\documentclass[12pt, leqno]{amsart}

\usepackage{graphicx}
\usepackage{amssymb,amsmath,amsthm,amscd,accents}
\usepackage{mathrsfs}
\usepackage{lscape}
\usepackage{enumerate}
\usepackage{upgreek}
\usepackage[usenames,dvipsnames]{color}
\usepackage[colorlinks=true, pdfstartview=FitV,
 linkcolor=blue,citecolor=blue,urlcolor=blue]{hyperref}
\usepackage{tikz}
\tikzset{dynkdot/.style={circle,draw,scale=.38}}
\usetikzlibrary{arrows.meta,arrows}
\usepackage{bm}
\usepackage{marginnote}
\usepackage{verbatim}

\usepackage[normalem]{ulem}  
\usepackage{xspace}
\usepackage{yfonts}
\usepackage{stmaryrd}

\usepackage[all]{xy}

\allowdisplaybreaks[3]

\newcommand{\arxiv}[1]{\href{http://arxiv.org/abs/#1}{\texttt{arXiv:#1}}}

\newcommand{\nc}{\newcommand}
\numberwithin{equation}{section}

\nc{\hs}{\hspace*}
\nc{\ms}{\mspace}

\nc{\qR}[1]{\ttq_{\mspace{-2mu}\raisebox{-.8ex}{${\scriptstyle{#1}}$}}}

\theoremstyle{plain}
\newtheorem{lemma}{Lemma}[section]
\newtheorem{proposition}[lemma]{Proposition}
\newtheorem{theorem}[lemma]{Theorem}

\newtheorem{corollary}[lemma]{Corollary}

\newtheorem*{convention}{Convention}

\theoremstyle{definition}
\newtheorem{remark}[lemma]{Remark}

\newtheorem{definition}[lemma]{Definition}

\newcommand{\longepito}[1][]{\xymatrix@C=4ex{{}\ar@{->>}[r]^{#1}&{}}}

\renewcommand{\le}{\leqslant}
\renewcommand{\ge}{\geqslant}
\renewcommand{\preceq}{\preccurlyeq}

\newcommand{\Ker}{{\operatorname{Ker}}}
\newcommand{\st}[1]{\{{#1}\}}

\newcommand{\seteq}{\mathbin{:=}}

\newcommand{\soplus}{\mathop{\mbox{\normalsize$\bigoplus$}}\limits}
\newcommand{\sbcup}{\mathop{\mbox{\normalsize$\bigcup$}}\limits}
\newcommand{\sbcap}{\mathop{\mbox{\normalsize$\bigcap$}}\limits}

\newcommand{\tens}{\mathop\otimes}

\newcommand{\g}{\mathfrak{g}}

\newcommand{\n}{\mathfrak{n}}

\newcommand{\C}{\mathbb{C}}
\newcommand{\Q}{\mathbb{Q}}
\newcommand{\Z}{\mathbb{Z}\ms{1mu}}

\newcommand{\al}{{\ms{1mu}\alpha}}
\newcommand{\ep}{\epsilon}
\newcommand{\la}{\lambda}
\newcommand{\be}{{\ms{1mu}\beta}}
\newcommand{\ga}{\gamma}

\newcommand{\La}{\Lambda}

\newcommand{\vph}{\varphi}

\newcommand{\wt}{{\rm wt}}
\newcommand{\hwt}{\widehat{\wt}}
\newcommand{\cl}{\colon}
\newcommand{\fin}{{\ms{1.5mu}{\rm fin}}}

\newcommand{\het}{{\rm ht}}

\newcommand{\Dynkin}{\triangle}  

\newcommand{\im}{\imath}
\newcommand{\jm}{\jmath}

\newcommand{\ii}{ \textbf{\textit{i}}}

\newcommand{\Hom}{\operatorname{Hom}}

\newcommand{\End}{\operatorname{End}}

\newcommand{\tY}{\widetilde{Y}}

\newcommand{\te}{\widetilde{e}}
\newcommand{\tf}{\widetilde{f}}

\newcommand{\hB}{\widehat{B}}

 \newcommand{\ocalD}{\overline{\calD}}

\newcommand{\fraks}{\mathfrak{s}}

\newcommand{\sfC}{\mathsf{C}}
\newcommand{\sfD}{\mathsf{D}}
\newcommand{\sfE}{\mathsf{E}}

\newcommand{\sfL}{\mathsf{L}}

\newcommand{\sfP}{\mathsf{P}}
\newcommand{\sfQ}{\mathsf{Q}}

\newcommand{\sfW}{\mathsf{W}}

\newcommand{\sfc}{\mathsf{c}}
\newcommand{\sfd}{\mathsf{d}}

\newcommand{\sfm}{\mathsf{m}}

\newcommand{\bbA}{\mathbb{A}}

\newcommand{\bbE}{\mathbb{E}}

\newcommand{\bbK}{\mathbb{K}}

\newcommand{\bbM}{\mathbb{M}}

\newcommand{\bzeta}{{\zeta_{\raisebox{-.2ex}{$\scriptstyle -$}}}}
\newcommand{\bzetai}[1][i]{{\zeta_{\raisebox{-.2ex}{$\scriptstyle{#1}, -$}}}}
\newcommand{\ubzetai}[2]{{\zeta_{\raisebox{-.2ex}{$\scriptstyle{#1}, -$}}^{#2}}}
\newcommand{\ubzeta}[1][-1]{{\zeta_{\raisebox{-.2ex}{$\scriptstyle -$}}^{#1}}}

\newcommand{\bfL}{\mathbf{L}}

\newcommand{\bfa}{\mathbf{a}}
\newcommand{\bfb}{\mathbf{b}}
\newcommand{\bfc}{\mathbf{c}}

\newcommand{\bfe}{\mathbf{e}}

\newcommand{\bfi}{\mathbf{i}}

\newcommand{\bfk}{\mathbf{k}}

\newcommand{\bfu}{\mathbf{u}}

\newcommand{\calA}{\mathcal{A}}

\newcommand{\calD}{\mathcal{D}}

\newcommand{\calK}{\mathcal{K}}

\newcommand{\calM}{\mathcal{M}}

\newcommand{\calQ}{\mathcal{Q}}

\newcommand{\calY}{\mathcal{Y}}

\newcommand{\hcalA}{\widehat{\calA}}

\newcommand{\scrB}{\mathscr{B}}
\newcommand{\scrC}{\mathscr{C}}
\newcommand{\scrD}{\mathscr{D}}

\newcommand{\scrU}{\mathscr{U}}

\newcommand{\ttP}{\mathtt{P}}

\newcommand{\ttb}{\mathtt{b}}

\newcommand{\ttq}{\mathtt{q}}

\newcommand{\rmE}{\mathrm{E}}

\newcommand{\rmG}{\mathrm{G}}

\newcommand{\rmL}{\mathrm{L}}

\newcommand{\rmP}{\mathrm{P}}

\newlength{\mylength}
\newcommand*{\para}{%
  \rlap{\rotatebox{-30}{\rule[.05ex]{.4pt}{.77em}}}%
  \kern.04em%
  \rlap{\kern.36em\raisebox{0.649519052835em}{\rule{.6em}{.4pt}}}%
  \rule{.6em}{.4pt}\kern-.04em%
  \rotatebox{-30}{\rule[.05ex]{.4pt}{.77em}}}

\newcommand{\isoto}[1][]{\mathop{\xrightarrow%
[{\raisebox{.3ex}[0ex][.3ex]{$\scriptstyle{#1}$}}]%
{{\raisebox{-.6ex}[0ex][-.6ex]{$\mspace{2mu}\sim\mspace{2mu}$}}}}}

\newcommand{\wl}{\sfP}
\newcommand{\rl}{\sfQ}
\newcommand{\weyl}{\sfW}
\newcommand{\lan}{\langle}
\newcommand{\ran}{\rangle}
\newcommand{\Bqg}{\scrB_q(\g)}
\newcommand{\Aqn}{A_q(\n)}
\newcommand{\Aqhn}{A_{q^{1/2}}(\n)}
\newcommand{\Azn}{A_\Zq(\n)}

\newcommand{\qt}[1]{\quad\text{#1}}
\newcommand{\qtq}[1][{and}]{\quad\text{{#1}}\quad}

\newcommand{\ee}{\end{enumerate}}
\newcommand{\bitem}{\begin{itemize}}
\newcommand{\eitem}{\end{itemize}}
\newcommand{\ben}{\begin{enumerate}[{\rm (1)}]}
\newcommand{\bnum}{\begin{enumerate}[{\rm (i)}]}
\newcommand{\bnump}{\begin{enumerate}[{\rm (i)$'$}]}
\newcommand{\bna}{\begin{enumerate}[{\rm (a)}]}
\newcommand{\bnA}{\begin{enumerate}[{\rm (A)}]}
\newcommand{\bc}{\begin{cases}}
\newcommand{\ec}{\end{cases}}
\newcommand{\ba}{\begin{array}}
\newcommand{\ea}{\end{array}}

\newcommand{\snoi}{\smallskip \noindent}
\newcommand{\mnoi}{\medskip \noindent}

\newenvironment{myequation}
{\relax\setlength{\arraycolsep}{1pt}\begin{eqnarray}}
{\end{eqnarray}}
\newenvironment{myequationn}
{\relax\setlength{\arraycolsep}{1pt}\begin{eqnarray*}}
{\end{eqnarray*}}

\nc{\eqs}[1]{\underset{\raisebox{.4ex}[.7ex][0ex]{$\scriptstyle{#1}$}}{=}}

\newcommand{\UZ}{\scrU_\Z}
\newcommand{\gf}{{\g_{\ms{1mu}{\rm fin}}}} 
\newcommand{\hBi}{\hB(\infty)}
\newcommand{\trmG}{\widetilde{\rmG}}
\newcommand{\LuphA}{\Lup\big(\hcalA_\Zq\big)}
\newcommand{\hBel}{\bfb = (b_k)_{k\in \Z}}
\newcommand{\Uqgm}{U_q^-(\g)}
\newcommand{\Uqhm}{U_{q^{1/2}}^-(\g)}
\newcommand{\Uzgm}{U_\Zq^-(\g)}
\newcommand{\Deltan}{\Delta_\n}

\newcommand{\es}{e^*}
\newcommand{\Es}{\rmE^\star}
\newcommand{\Esn}[1]{\rmE^{\star \hspace{0.1ex} #1}}
\newcommand{\Zq}{{\Z[q^{\pm1}]}}
\newcommand{\Zqh}{{\Z[q^{\pm1/2}]}}
\newcommand{\Qh}{{\Q(q^{1/2})}}
\newcommand{\Qq}{{\Q(q)}}

\newcommand{\Gup}{\rmG^{{\rm up}}}
\newcommand{\Lup}{\rmL^{{\rm up}}}
\newcommand{\LupA}{\Lup\left(\Azn\right)}
\newcommand{\rdual}{\scrD}

\newcommand{\llf}{(\hspace{-0.5ex}(}

\newcommand{\rrf}{)\hspace{-0.5ex})}

\newcommand{\hAm}[1]{  \hcalA[#1]  }
\newcommand{\hAmz}[1]{  \hcalA[#1]_\Zq  }
\newcommand{\Ang}[1]{ \bigl\langle #1 \bigr\rangle}
\newcommand{\dVert}[1]{  \Vert #1 \Vert  }
\newcommand{\hAform}[1]{\bl #1\br_{\hA}\ms{1mu}}
\newcommand{\kform}[1]{  ( #1 )}
\newcommand{\aform}[1]{  ({ #1 })_\n  }
\newcommand{\Aform}[1]{  \boldsymbol{\left(} #1 \boldsymbol{\right)}_\n  }
\newcommand{\bAng}{\Ang}

\newcommand{\pair}[1]{  \llf #1 \rrf  }

\newcommand{\Iset}[1]{  \{ #1 \}_{i \in I} }
\newcommand{\Cg}{\scrC_\g}

\newcommand{\akew}[1][2ex]{\rule[-1ex]{#1}{0ex}}
\newcommand{\oprod}{\prod\limits^{{\scriptstyle\xrightarrow{\akew[1ex]}}}}
\newcommand{\rprod}{\prod^{\xleftarrow{}}}
\newcommand{\Qdatum}{(\Dynkin,\varsigma,\xi)} 
\newcommand{\Mn}{\mathbf{M}\ms{1mu}}

\newcommand{\xc}{{x_\circ}}
\newcommand{\yc}{{y_\circ}}

\newcommand{\fc}{{f_\circ}}

\newcommand{\akete}[1][0ex]{\rule[{#1}]{0ex}{1ex}}

\newcommand{\sumI}{\sum_{i\in I}}

\newcommand{\eq}{\begin{myequation}}
\newcommand{\eneq}{\end{myequation}}
\newcommand{\eqn}{\begin{myequationn}}
\newcommand{\eneqn}{\end{myequationn}}
\newcommand{\bnom}[1]{\begin{bmatrix}#1\end{bmatrix}}
\newcommand{\vs}{\vspace*}

\newcommand{\bl}{\bigl(}
\newcommand{\br}{\bigr)}

\newcommand{\epito}{\twoheadrightarrow}

\nc{\hAg}[1][m]{\hcalA_{\ge{#1}}} 
\nc{\hAsg}[1][m]{\hcalA_{>{#1}}} 
\nc{\hAl}[1][m]{\hcalA_{\le{#1}}} 
\nc{\hAsl}[1][m]{\hcalA_{<{#1}}}
\newcommand{\hA}{\hcalA}
\newcommand{\Um}{\Uqgm}

\newcommand{\vphi}{\vph}
\newcommand{\ang}[1]{\langle#1\rangle}
\nc{\catC}{{\mathscr{C}}} 		
\nc{\catCQ}{\catC_{\cQ}}
\nc{\catCO}{\catC_{\g}^0}
\nc{\Zp}{\Z_{\ge0}}
\nc{\Zsp}{\Z_{>0}} 
\nc{\tin}{t^{-1}}  
\nc{\En}[1][i]{{\mathrm{E}^\n_{\ms{2mu}{#1}}}}
\nc{\Ens}[1][i]{{\mathrm{E}^{*\,\n}_{\ms{2mu}{#1}}}}
\nc{\qq}{\mathsf{q}}
\nc{\Umz}{\Um?}
\nc{\wtUm}{\Um??}
\nc{\nrl}{\rl^-}
\nc{\ol}{\overline}
\nc{\hcalAf}{\hcalA}
\nc{\scbul}{{\scriptstyle\bullet}}

\usepackage[colorinlistoftodos]{todonotes}

\title[Global bases for Bosonic extensions]{Global bases for Bosonic extensions of quantum unipotent coordinate rings}

\author[M. Kashiwara]{Masaki Kashiwara}
\thanks{The research of M.\ Kashiwara
	was supported by Grant-in-Aid for Scientific Research (B)  23K20206,  
	Japan Society for the Promotion of Science.}
\address[M. Kashiwara]{%
Kyoto University Institute for Advanced Study, Research Institute
for Mathematical Sciences, Kyoto University, Kyoto 606-8502, Japan
\& Korea Institute for Advanced Study, Seoul 02455, Korea }
\email[M. Kashiwara]{masaki@kurims.kyoto-u.ac.jp}

\author[M. Kim]{Myungho Kim}
\address[M. Kim]{Department of Mathematics, Kyung Hee University, Seoul 02447, Korea}
\email[M. Kim]{mkim@khu.ac.kr}
\thanks{The research of M.\ Kim was supported by the National Research Foundation of
Korea (NRF) Grant funded by the Korea government(MSIT)
NRF-2022R1F1A1076214 and NRF-2020R1A5A1016126.}

\author[S.-j. Oh]{Se-jin Oh}
\thanks{ The research of S.-j.\ Oh was supported by the National Research Foundation of
	Korea (NRF) Grant funded by the Korea government(MSIT) (NRF-2022R1A2C1004045).}
\address[S.-j. Oh]{ Department of Mathematics, Sungkyunkwan University, Suwon, South Korea}
\email[S.-j. Oh]{sejin092@gmail.com}

\author[E. Park]{Euiyong Park}
\thanks{The research of E.\ Park was supported by the National Research Foundation of Korea (NRF) Grant funded by the Korea Government(MSIT)(RS-2023-00273425 and NRF-2020R1A5A1016126).}
\address[E. Park]{Department of Mathematics, University of Seoul, Seoul 02504, Korea}
\email[E. Park]{epark@uos.ac.kr}

\keywords{Bosonic extension, Global basis, Hernandez-Leclerc category, Quantum Grothendieck ring} 

\subjclass[2010]{05E10, 05E18, 17B37} %

\date{2024, June 18}

\begin{document}

\begin{abstract}
In the paper, we establish the global basis theory for the bosonic extension $\hA$ associated with an arbitrary generalized Cartan matrix.
When $\hA$ is of simply-laced finite type, $\hA$ is isomorphic to the quantum Grothendieck ring $\calK_q(\Cg^0) $ of the Hernandez-Leclerc category $\catCO$ over a quantum affine algebra.
In this case, we show that the $(t,q)$-characters  of simple modules in $ \catCO$ correspond to the normalized global basis of $\hA$.	
\end{abstract}

\setcounter{tocdepth}{1}

\maketitle
\tableofcontents

\section*{Introduction}

Let $U_q'(\g)$ be a quantum affine algebra of untwisted type, where $q$ is an indeterminate, and let $\catC_\g$ be  the category of finite-dimensional integrable $U_q'(\g)$-modules. 
The category $\catC_\g$ has a rich and beautiful structure including rigidity and has been actively studied in various aspects (see \cite{CP95,FR99, HL10,K02,
KKOP22, KKOP24, Nak04, VV03}  and references therein).  
There exists a distinguished monoidal subcategory $\catCO$ of $\catC_\g$, called a \emph{Hernandez-Leclerc category}, which is  defined by certain fundamental modules (\cite{HL10}). The Hernandez-Leclerc category $\catCO$ takes a central position in studying $\catC_\g$.
 Indeed, the Grothendieck ring $K(\catC_\g)$ is isomorphic to the tensor products of
infinitely many copies of $K(\catC^0_\g)$. 
It was proved that various subalgebras of the Grothendieck rings $K(\catCO)$ have cluster algebra structures (\cite{HL10, HL15, FHOO23}) and the corresponding subcategories provide monoidal categorifications (\cite{KKOP24}).

The \emph{$q$-characters} of simple modules in $\catC_\g$ were introduced by Frenkel and Reshetikhin as an analogue of the ordinary characters for finite type quantum groups (\cite{FR99}). 
For a simple module $L$ in $\catCO$, we denote by  $\chi_{q}(L)$ the $q$-character of $L$. 
The $q$-characters lead us to understand the Grothendieck ring $K(\catCO)$ as a subalgebra of the Laurent polynomial ring $\calY = \Z[Y_{\pi(\im),p}^{\pm 1} \mid (\im,p) \in \sigma_0(\g)]$ (see Section \ref{Sec: qGr} for the definition of $\sigma_0(\g)$). The Grothendieck ring $K(\catCO)$  is the  commutative polynomial ring generated by the isomorphism classes of \emph{fundamental modules} in $\catCO$.

The $t$-analogues of $q$-characters of simple modules, called \emph{$(q,t)$-characters}, were introduced in \cite{Nak04, VV03} for simply-laced types and in \cite{Her04} for non-simply laced types. 
For a simple module $L$ in $\catCO$, we denote by   $\chi_{q,t}(L)$ the $(q,t)$-character of $L$. 
The $(q,t)$-characters $\chi_{q,t}(L)$ allow us to consider the \emph{quantum Grothendieck ring} $\calK_t(\Cg^0)$ as a subalgebra of the \emph{quantum torus}  $\calY_t$ (see Section  \ref{Sec: qGr} for the definition of $\calY_t$). 
The quantum Grothendieck ring $\calK_t(\Cg^0)$ is non-commutative. 
When specializing at $t=1$, the $(q,t)$-characters $\chi_{q,t}(L)$ correspond to the $q$-characters $\chi_q(L)$ for a simple module $L \in \catCO$ satisfying either that
$\g$ is of untwisted affine  $ADE$ type (\cite{Nak04, VV03}) and  $B_n^{(1)}$ type (\cite{FHOO22})  
or that $\g$ is of general type but $L$ is \emph{reachable} (\cite{FHOO23}).
The set of all $(q,t)$-characters $\chi_{q,t}(L)$ of simple modules $L \in \catCO$ becomes a natural $\Z[t^{\pm 1/2}]$-linear basis of $\calK_t(\Cg^0)$ with several positivity properties. This basis is sometimes  called a \emph{canonical basis} of $\calK_t(\Cg^0)$ (see \cite{FHOO22, FHOO23} for example), 
and a connection to the \emph{upper global basis} (or \emph{dual canonical basis}) of $U_q^-(\gf)$ (see \cite{Kashiwarabook, LusztigBook} and references therein) is shown in \cite{FHOO22, HL10} via the subcategories $\catC_\calQ$.
On the combinatorial side, the set of isomorphism classes of simple modules in $\catCO$ has a \emph{categorical crystal structure} and it is isomorphic to the \emph{extended crystal} $\hB(\infty)$ of the infinite crystal $B(\infty)$ of $U_q^-(\gf)$, where $\gf$ is given in Section \ref{Sec: qGr}  (see \cite{KP22, Park23} for categorical and extended crystals).

A presentation of the quantum Grothendieck ring $\calK_t(\Cg^0)$
as a ring is given in \cite{HL15} (see also \cite{FHOO22}). 
The ring $\calK_t(\Cg^0)$ has an infinitely many generators  $ \{ f_{i,p} \mid i\in I_\fin, p\in \Z \}$ and the \emph{quantum Serre relations} and the \emph{quantum bosonic relations} hold for the generators $f_{i,p}$ ($i\in I_\fin$, $p\in \Z$). For each $p\in \Z$, the subalgebra of  $\calK_t(\Cg^0)$ generated by $f_{i,p}$ ($i\in I_\fin$) is isomorphic to the negative half $U_t^-(\gf)$ or more precisely the \emph{quantum unipotent coordinate ring} $A_t(\n)$. Thus the whole algebra  $\calK_t(\Cg^0)$ can be understood as an \emph{affinization} of $A_t(\n)$.
When $\g$ is of  untwisted affine $ADE$ type,
he quantum Grothendieck ring $\calK_t(\Cg^0)$ has a connection to the derived Hall algebra of the category of representations of a quiver $Q$ (see \cite{HL15}).

Given a \emph{finite type Cartan matrix} $\sfC$,  the \emph{bosonic extension} $\hA$ associated with $\sfC$  is the algebra formally defined  by the generators $f_{i,p}$ ($i\in I$, $p\in \Z$) with  the quantum Serre relations and quantum bosonic relations determined by the Cartan matrix $\sfC$. 
By construction,  the quantum Grothendieck ring $\calK_t(\Cg^0)$ is isomorphic to the bosonic extension $\hA$ associated with $\gf$. Note that $\gf$
is always of simply-laced finite type. When $\sfC$ is of not simply-laced finite  type, $\hA$ is isomorphic to a \emph{quantum virtual Grothendieck ring} realized inside a quantum torus (\cite{JLO23A, JLO23B}). In \cite{OP24}, the third and forth named authors introduced a new family of subalgebras $\hA(\ttb)$ of $\hA$ for any elements $\ttb$ in the positive submonoid of the (generalized) \emph{Braid group} associated with $\sfC$. They developed the  \emph{Poincar\'e-Birkhoff-Witt theory} for $\hA (\ttb)$ using the \emph{braid group actions} $\mathsf{T}_i$ on $\hA$ (see \cite{KKOP21A} and \cite{JLO23B} for braid group actions). It is conjectured that the algebra $\hA(\ttb)$ has a quantum cluster algebra structure.

In this paper, we establish the \emph{global basis theory} for the bosonic extension $\hA$ associated with an \emph{arbitrary generalized Cartan matrix}. 
The main results can be summarized briefly as follows.
\bnum
\item We introduce the bosonic extension $\hA$ associated with an \emph{arbitrary generalized Cartan matrix} (see \eqref{Eq: def of hA}). We prove that there exists a \emph{faithful} representation of $\hA$ (Theorem \ref{Thm: repn of hA}), which allows us to realize $\hA$.
  Note that  a quantum torus is not used to realize $\hA$. 
  In the definition of $\hA$, we use $q$ instead of $t$ 
  as an indeterminate for consistency of the previous works by the authors. 
 As a consequence, we have the serial decomposition of $\hA$ (Corollary \ref{cor: main1}), which
allows us to understand $\hA$  as an \emph{affinization} of the quantum unipotent coordinate ring $A_q(\n)$ associated with $\sfC$.

\item We construct the \emph{global basis} of $\hA$. We first construct the $\Zq$-lattice $\hcalA_\Zq$ of $\hA$ by using the lattice $\Azn$  generated by the \emph{upper global basis} of $A_q(\n)$. We then prove that there exists a distinguished linear  basis of  $\hcalA_\Zq$  with the same properties as the upper global basis of  $A_q(\n)$ (Theorem \ref{Thm: global basis}). This distinguished basis is parameterized by the extended crystal $\hB(\infty)$ and is orthogonal with respect to the \emph{symmetric bilinear form} $ \hAform{ \  , \   }$ on $\hA$ at $q=0$  (see \eqref{eq: hA form},  Definition \ref{def: global basis} and Proposition~\ref{Prop: properties of Gb}).

\item When $ \hA$ is of simply-laced finite type, there is an isomorphism between $ \hcalAf $  and  $ \calK_q(\Cg^0) $.
Under this isomorphism,  we show that the $(t,q)$-characters $\chi_{t,q}(L)$ of simple modules $L \in \catCO$ correspond to the \emph{normalized global basis} of $\hA$  (see Theorem \ref{Thm: simple to gb}). Here,  we swap the role of $t$ and $q$, and use $q$ in $\calK_q(\Cg^0)$ for consistency (see Section \ref{Sec: qGr}). 
\ee

Let us explain the main results more precisely. Let $\sfC=(c_{i,j})_{i,j\in I}$ be a generalized Cartan matrix of \emph{symmetrizable Kac-Moody type}, and let $U_q(\g)$ (resp.\ $A_q(\n)$) be the quantum group (resp.\ quantum unipotent coordinate ring) associated with $\sfC$. We define $\hcalA$ to be the $\Qh$-algebra generated by $\{f_{i,p}\mid i\in I,p\in \Z\}$ with the defining relations \eqref{Eq: def of hA}.
For $ - \infty \le a \le b \le \infty$, we denote by $\hcalA[a,b]$ the subalgebra of $\hcalA$ generated by $\{ f_{i,k} \ | \  i \in I, a\le k \le b\}$, and set 
$\hcalA[m] \seteq \hcalA[m,m]$.
Define
$$
\UZ = \sbcup_{a\le b}\UZ[a,b],
$$
where $\scrU_\Z[a,b]=\Uqhm_b\tens\Uqhm_{b-1}\tens\cdots\tens\Uqhm_{a+1}\tens\Uqhm_a.$
Theorem \ref{Thm: repn of hA} tells us that the action of $f_{i,m}$ given in \eqref{eq: action f_im} gives a faithful  $\hcalA$-module structure on $\scrU_\Z$. In the course of the  proof of Theorem \ref{Thm: repn of hA} , the quantum Boson algebra $\Bqg$ plays a crucial role (see \eqref{eq: Phi bzeta} and Lemma \ref{lem: Ezeta}). As a consequence, Corollary \ref{cor: main1} says that $\hcalA[ m]$ is isomorphic to $\Aqhn \simeq \Uqhm  $ for each $m\in \Z$ and that $\hA[a,b]$ is decomposed as follows:
$$
\hcalA[a,b] \simeq  \hcalA[b] \tens_\Qh \hcalA[b-1] \tens_\Qh \cdots \tens_\Qh \hcalA[a+1] \tens_\Qh \hcalA[a].
$$

This decomposition provides us the natural projection
$	\Mn\colon \hcalA\longepito \oprod_{k\in\Z}\hcalA[k]_{0}\simeq\Qh$.
We define a bilinear form on $\hcalA$ by
\begin{align*}
	\hAform{x,y} \seteq \Mn(x \ocalD(y)) \in \Qh \quad \text{ for any } x,y \in \hcalA,   
\end{align*}
where $\ocalD$ is the automorphism of $\hcalA$ given in~\eqref{Eq: shift}.
We introduce the operators $\rmE_{i,m}$ and $\Es_{i,m}$ on $\hA$ using the $q$-brackets (see \eqref{eq: Ei Esi}) and investigate their properties. 
Lemma \ref{lem: E Es} says that $\rmE_{i,m}$ and $\Es_{i,m}$ play the role of $e_i'$ and $e_i^*$ in the local piece $\hA[m]$. 
We then prove that $\hAform{\ ,\ }$ is symmetric and non-degenerate  (see Theorem \ref{thm: hAform}). 
Interestingly, $f_{i,m}$ and $f_{i,m+1}$ acts as adjoint operators, i.e., 
$$
\hAform{f_{i,m}x, y } =  \hAform{x, y f_{i,m+1} }
$$ 
for any  $x,y\in \hA$ (see Lemma \ref{lem: hAform}). The bilinear form $\hAform{\ , \ }$ induces another bilinear form $\pair{\ , \ }$ on $\hA$ by twisting the power of $q$ (see \eqref{Eq: def of pair}).
This bilinear form $\pair{\ , \ }$ coincides with the bilinear form on $\hA$ introduced in \cite{OP24} when $\sfC$ is of finite type.

For each $k \in \Z$, there exists a $\Qh$-algebra isomorphism 
$$
\varphi_{k} \cl \Aqhn   \isoto  \hcalA[k]
$$
(see \eqref{eq: vph} for the definition of $\varphi_{k}$). Let $\sfc$ be the endomorphism of $\Aqhn$ defined in \eqref{eq: sfc}, which is used in the characterization of the upper global basis of $A_q(\n)$.
The bosonic extension $\hA$ has two kind of maps  $\overline{\phantom{a}} \in \End(\hA)$ and $c \in \End(\hA)$ (see Section \ref{Sec: BE}), which are used to define the (normalized) global basis of $\hA$. Note that they only differ up to a power of $q$ (see \eqref{eq: c map}) 

We first show that the isomorphism $\varphi_{k}$ is compatible with the map $\sfc \in  \End(\Aqhn) $ and the map $c \in \End(\hA)$ (see Lemma \ref{lem: vph_k}).
Define
$$
\hcalA[k]_\Zq \seteq \vph_{k}(\Azn) \subset \hcalA,
$$
where $\Azn$ is the $\Zq$-lattice generated by the upper global basis of $A_q(\n)$, and set
$$
\hcalA[a,b]_\Zq \seteq \oprod_{k \in [a,b]} \hcalA[k]_\Zq\subset \hcalA,  \qquad 
\hcalA_\Zq \seteq \sbcup_{a \le b} \hcalA[a,b]_{\Zq}\subset \hcalA. 
$$
Proposition \ref{prop: hcalA integral form} says that $\hcalA_\Zq$ is a subalgebra invariant under the action of $c$.
We denote by $\hB(\infty)$ the extended crystal of the infinite crystal $B(\infty)$ of $U_q^-(\g)$ (see Section \ref{subsec: crystals} for details).
For any $\bfb =(b_k)_{k\in\Z} \in \hB(\infty)$, we set
$$
\rmP(\bfb) \seteq \oprod_{k \in \Z} \vph_k(\Gup(b_k)) \in \Lup(\hcalA_\Zq), 
$$
where $\Lup\left( \hcalA_\Zq \right)$ is the $\Z[q]$-lattices inside $\hcalA_\Zq$ defined in \eqref{Eq: def of Zq lattices of hA}, and $\Gup(b_k)$ denotes the member of the upper global basis of $A_q(\n)$ corresponding to $b_k$. We then prove the following.
\begin{theorem}  [= Theorem \ref{Thm: global basis}] \ 
\bnum
\item 
For each $\hBel \in\hBi$, there exists a unique
$\rmG(\bfb)\in \LuphA$ such that
\bna
\item
  $\rmG(\bfb)-\rmP(\bfb)\in \displaystyle\sum_{\bfb'\prec\bfb}q\Z[q]\rmP(\bfb')$, where $\prec$ is a certain order on $\hBi$ {\rm (}see \eqref{Eq: order on hB}{\rm)}.
\item
$c(\rmG(\bfb))=\rmG(\bfb)$.  
\ee
\item Moreover, the set $\{\rmG(\bfb) \ | \ \bfb\in\hBi\}$ forms a $\Z[q]$-basis of $\LuphA$, and
a $\Z$-basis of $\LuphA\cap c\big(\LuphA\big)$. 
\item 
For any $\bfb \in\hBi$, we have 
$$
\rmP(\bfb) = \rmG(\bfb) + \sum_{\bfb'  \prec \bfb}  f_{\bfb,\bfb'}(q) \rmG(\bfb')\qquad \text{ for some $f_{\bfb,\bfb'}(q) \in q \Z[q]$.}
$$
\ee
\end{theorem}

We define
$$
\rmG :=  \{\rmG(\bfb)   \mid   \bfb\in\hBi \} \quad \text{ and } \quad \trmG:=  \{ \trmG(\bfb)    \mid   \bfb\in\hBi \}
$$ 
where $ \trmG(\bfb) \seteq q^{ - (\wt( \rmG(\bfb),\,\wt( \rmG(\bfb)  )/4} \rmG(\bfb) $.
 We call  $\rmG$  and $ \trmG$  the \emph{global basis}  and the \emph{normalized global basis} of $\hcalA$ respectively.
Note that  $\trmG(\bfb)$ is bar-invariant by \eqref{Eq: bar inv c inv}.
It turns out that the (normalized) global basis of $\hA$ enjoys the same properties as the upper global basis of $A_q(\n)$.
In particular, the normalized global basis is orthogonal with respect to $ \hAform{ \  , \   }$ at $q=0$ (see Proposition \ref{Prop: properties of Gb}).

We next study the global basis in the viewpoint of the quantum Grothendieck ring. 
Let $U_t'(\g)$ be the quantum affine algebra of untwisted affine  type and let $\catCO$ be its Hernandez-Leclerc category. 
We choose a  $Q$-datum $\calQ$ and denote by $\hA$ the bosonic extension associated with $\gf$. Then there is
 an $\Qh$-algebra isomorphism
\begin{align*}
	\Omega_\calQ \cl  \bbK^0  \buildrel \sim \over \longrightarrow \hcalAf,
\end{align*}	
where $\bbK^0 \seteq   \Qh \otimes_{\Zqh } \calK_q(\Cg^0)$ (see \eqref{Eq: OmegaQ}  for the definition of $\Omega_\calQ$).
We identify the extended crystal $\hB (\infty)$ with the set $\Z^{\oplus \Z}_{\ge0}$ via the \emph{PBW monomial realization} associated with a $\calQ$-adapted reduced sequence $\bfi$ of  the longest element $w_0 \in \weyl_\fin$ (see Section \ref{subsec: crystals} for details). Let $\sfm(\bfc)$ be the \emph{dominant monomial} corresponding to $\bfc \in \Z^{\oplus \Z}_{\ge0}$ and let $ L( \bfc)$ be the simple module in $\catCO$ corresponding to  the dominant monomial $\sfm(\bfc)$ (see Section \ref{Sec: gb sm} for details). Utilizing the unitriangularity, we prove that the $ (t,q)$-characters $\chi_{t,q}( L( \bfc) )$ correspond to $\trmG(\bfc)$ under the isomorphism $\Omega_\calQ$ (see Theorem \ref{Thm: simple to gb}).

\medskip
This paper is organized as follows. 
Section \ref{Sec: Preliminaries}  is devoted to give necessary backgrounds on quantum groups and crystals. 
In Section \ref{Sec: QBA} and Section \ref{Sec: ugb}, we investigate several properties of the quantum bosonic algebras and the upper global bases. 
In Section \ref{Sec: BE}, we introduce and realize the bosonic extension $\hA$ associated with an arbitrary generalized Cartan matrix.
In Section \ref{Sec: hAform}, we construct and study the bilinear form $\hAform{ \ ,\ }$ on $\hA$ 
and, in Section \ref{Sec: GB}, we construct the global basis of $\hA$. 
In Section  \ref{Sec: QGR and GB}, we finally show that the $(t,q)$-characters correspond to the normalized global basis.

\medskip

\vskip 1em

{\bf Acknowledgments}
The second, third and fourth authors gratefully acknowledge for the hospitality of RIMS (Kyoto University) during their visit in 2024.

\vskip 2em

\begin{convention}  Throughout this paper, we use the following convention.
\ben
\item For a statement $\ttP$, we set $\delta(\ttP)$ to be $1$ or $0$ depending on whether $\ttP$ is true or not. In particular, we set $\delta_{i,j}=\delta(i=j)$. 
\item For a totally ordered set $J = \{ \cdots < j_{-1} < j_0 < j_1 < j_2 < \cdots \}$, write
$$
\oprod_{j \in J} A_j \seteq \cdots A_{j_2}A_{j_1}A_{j_0}A_{j_{-1}}A_{j_{-2}} \cdots, \quad
\rprod_{j \in J} A_j \seteq \cdots A_{j_{-2}}A_{j_{-1}}A_{j_0}A_{j_{1}}A_{j_{2}} \cdots. 
$$
\item For $a\in \Z \cup \{ -\infty \} $ and $b\in \Z \cup \{ \infty \} $ with $a\le b$, we set 
\begin{align*}
	& [a,b] =\{  k \in \Z \ | \ a \le k \le b\}, &&  [a,b) =\{  k \in \Z \ | \ a \le k < b\}, \allowdisplaybreaks\\
	& (a,b] =\{  k \in \Z \ | \ a < k \le b\}, &&  (a,b) =\{  k \in \Z \ | \ a < k < b\},
\end{align*}
and call them \emph{intervals}. 
When $a> b$, we understand them as empty sets. For simplicity, when $a=b$, we write $[a]$ for $[a,b]$. For an interval $[a,b]$, we set
$A^{[a,b]}$ to be the product of copies of a set $A$ indexed by $[a,b]$, and 
$$
\Zp^{\oplus [a,b]} \seteq \{ (c_a,\ldots,c_b) \ | \ c_k \in \Z_{\ge 0} \text{ and  $c_k=0$ except for finitely many $k$'s} \}. 
$$
We define $A^{[a,b)}$,  $\Zp^{\oplus [a,b)}$, ..., etc.\ in a similar way. 
\item 
  For a pair of elements $x,y$ in a $\Zq$-module, we write
  $x\equiv_qy$ if $x=q^my$ for some $m\in\Z$.
  \label{it:equivq}
\item
  For a preorder $\preceq$ on a set $A$,
  we write $x\prec y$ if $x\preceq y$ holds but
  $y\preceq x$ does not hold.
  \label{conv:preorder}
\ee
\end{convention}

\vskip 2em 

\section{Preliminaries} \label{Sec: Preliminaries}
In this section, we review the basic definitions and properties of
quantum groups and quantum unipotent coordinate rings. We refer the readers to \cite{Kashiwarabook,LusztigBook} for more details. 

\subsection{Quantum groups} Let $(\sfC,\wl,\Pi,\wl^\vee,\Pi^\vee)$ be a symmetrizable \emph{Cartan datum} consisting of a generalized \emph{Cartan matrix} $\sfC=(c_{i,j})_{i,j \in I}$
 indexed by  an index set $I$, a \emph{weight lattice} $\wl$, a set of \emph{simple roots} $\Pi =\{\al_i \}_{i \in I} \subset P$, 
a \emph{coweight lattice} $\wl^\vee \seteq \Hom_{\Z}(\wl,\Z)$ and a set of
\emph{simple coroots} $\{ h_i \}_{i \in I} \subset \wl^\vee$. The datum satisfies
$\Ang{h_i,\al_j}=c_{i,j}$ for all $i,j\in I$.  We choose $\Iset{\La_i}$ such that $\Ang{h_j,\La_i}=\delta_{i,j}$ for $i,j \in I$ and call them
the \emph{fundamental weights}.

Note that there is a  $\Q$-valued  symmetric bilinear form $( \ , \ )$ on $\wl$ such that
$$
(\al_i,\al_i) \in 2\Zp \qtq \Ang{h_i,\la} = \dfrac{2(\al_i,\la)}{(\al_i,\al_i)} \quad \text{ for any } \la \in \wl.
$$
There exists a diagonal matrix
$\sfD = {\rm diag}( \sfd_i \seteq (\al_i,\al_i)/2 \ | \ i \in I )$ such that $\sfD \sfC$ is symmetric. 

We denote by $\rl = \bigoplus_{i \in I}\Z \al_i$. 
We set $\sfQ^+ \seteq \sum_{ i \in I } \Zp\al_i$ the \emph{positive root lattice} and
$\sfQ^- \seteq -\sfQ^+$ the \emph{negative root lattice}. For any
$\be = \sumI a_i\al_i \in \rl$, we set 
\begin{align} \label{eq: het dVert}
\het(\be) = \sumI |a_i| \in  \Zp   \qtq  \Vert  \beta \Vert =  \sumI |a_i|\al_i \in \rl^+.    
\end{align}

Recall that there exists a partial order $\le$ on $\rl$ defined by 
\begin{align} \label{eq: < order on Q}
\text{$\al \le \be$ if $\be -\al \in \rl^+ $ for $\al,\be \in \rl$. }
\end{align}

\smallskip

Let $q$ be an indeterminate with a formal square root $q^{1/2}$. We define
$$
[n]_{q} \seteq \dfrac{q^n-q^{-n} }{q-q^{-1}}, \quad   [n]_{q} ! \seteq \prod_{k=1}^n [k]_{q}
\qtq \begin{bmatrix} m \\ n     \end{bmatrix}_{q} \hspace{-1ex} \seteq
\dfrac{[m]_{q}!}{[n]_{q}!\,[m-n]_{q}!}\qtq q_i \seteq q^{\sfd_i}
$$
for $i \in I$ and $m \ge n \in \Zp$.
We simply write $[n]_i$ (resp.\ $[n]$) instead of $[n]_{q_i}$ (resp.\ $[n]_q$). 

We denote by $U_q(\g)$ the \emph{quantum group} over $\Q(q)$ which is associated with a Cartan datum $(\sfC,\wl,\Pi,\wl^\vee,\Pi^\vee)$ and
generated by Chevalley generators $e_i,f_i$ $(i\in I)$ and $q^{h}$ $(h \in \wl^\vee)$. Let $U_q^-(\g)$ be the subalgebra of
$U_q(\g)$ generated by $f_i$ $(i \in I)$ and $U_{\Zq}^-(\g)$ the $\Zq$-subalgebra of
$U_q(\g)$ generated by $f_i^{(n)} \seteq f_i^n/[n]_i!$ $(i \in I, n\in \Zsp)$. 

Remark that 
\bnum
\item $\Uqgm$ admits the \emph{weight space decomposition}
$$\Uqgm = \soplus_{\be \in  \rl^-} \Uqgm_\be.$$ 
For an element $x \in \Uqgm_\be$, we set $$\wt(x)= \be \in \rl^-.$$
\item There exists a $\Q(q)$-algebra anti-involution $*\cl U_q(\g) \to U_q(\g)$ defined by 
$$
e_i^* =e_i, \quad f_i^* =f_i \qtq (q^h)^* = q^{-h},
$$
where we write $x^*$ for the image of $x \in U_q(\g)$ under the involution $*$.
\ee

\subsection{Crystals} \label{subsec: crystals}
In this subsection, we briefly recall the notions of \emph{infinite crystals} $B(\infty)$ and \emph{extended crystals} $\hB(\infty)$.  
We refer the reader to \cite{K91, K93, K95, Kashiwarabook} for  the
details on crystals, and \cite{KP22, Park23} for the details on extended crystals.

Let $B(\infty)$ be the \emph{infinite crystal} of the negative half $\Uqgm$, and let $\tf_i$ and $\te_i$ be the \emph{crystal operators} for $B(\infty)$.   
For any $b\in B(\infty)$, $\wt(b)$ stands for the weight of $b \in B(\infty) $. 
The $\Qq$-algebra anti-involution $*\colon U_q(\g) \to  U_q(\g)$
induces the involution $*$ of $B(\infty)$
and defines another pair of crystal operators $\tf_i^*$ and $\te_i^*$ on $B(\infty)$:
$\tf_i^*b=\bl\tf_i(b^*)\br^*$ and $\te_i^*b=\bl\te_i(b^*)\br^*$.

The \emph{extended crystal} of $B(\infty)$ is defined as 
\begin{align} \label{Eq: extended crystal}
	\hB(\infty) \seteq   \Bigl\{  (b_k)_{k\in \Z } \in \prod_{k\in \Z} B(\infty)  \bigm | b_k =\mathsf{1} \text{ for all but finitely many $k$}  \Bigr\},
\end{align}
where $\mathsf{1}$ is the highest weight element of $B(\infty)$. The \emph{extended crystal operators} on $\hB (\infty)$ are defined by the usual crystal operators $\tf_i$, $\te_i$, $\tf_i^*$ and $\te_i^*$ (see \cite{KP22, Park23} for details).  

The involution 
${}^\star\cl \hB(\infty) \longrightarrow \hB(\infty)$ is defined as follows: for any $\bfb = (b_k)_{k\in \Z} \in \hB(\infty)$,   
\eq
\bfb^\star = (b_k')_{k\in \Z}\qt{where $b_k' := b_{-k}^*$   for any $ k\in \Z$ (see \cite[Section 4]{KP22}).}
\label{def:starc}
\eneq

\smallskip

When $\g$ is of finite type, one can realize $B(\infty)$ in a combinatorial way using \emph{PBW theory} on $\Uqgm$ (see \cite{SST18} and references therein). Let $\g$ be of finite type and let $\underline{w}_0 = s_{i_1} s_{i_2} \cdots s_{i_\ell} $ be a \emph{reduced expression} of the longest element $w_0$ in the Weyl group $\weyl$. We set 
$\ii \seteq (i_1, \ldots, i_\ell) \in I^\ell$, which is called a \emph{reduced sequence} of $w_0$.
 Matching $b \in B(\infty)$ with the corresponding \emph{PBW monomial} associated with $\ii$,
the infinite crystal $B(\infty)$ can be identified with the set $\Zp^{\oplus [1, \ell]}$, i.e., 
$$
B(\infty) \simeq \Zp^{\oplus [1,\ell]}.
$$ 
This PBW realization can be extended  to the extended crystal $\hB(\infty)$. For any $\bfc = (c_k)_{k\in \Z} \in \Zp^{\oplus \Z}$, and we  understand $\bfc_k \seteq(c_{k\ell+1}, c_{k\ell+2}, \ldots, c_{(k+1)\ell})$ as an element of $B(\infty)$ via the PBW realization. Considering the definition \eqref{Eq: extended crystal}, one can understand $\bfc$ as an element of $\hB(\infty) $. Thus we realize $\hB(\infty)$ using the PBW realization of $B(\infty)$, i.e., 
$$
\hB(\infty) \simeq \Zp^{\oplus \Z}.
$$
\subsection{Quantum unipotent coordinate rings}
In this subsection, we briefly review the definition of the \emph{quantum unipotent coordinate ring}. We refer the readers to \cite[Section 8]{KKKO18} for details.

Let us endow $\Uqgm \tens \Uqgm$ the algebra structure defined by 
$$
(x \tens y)(x' \tens y') = q^{-(\wt(x'),\wt(y))} (xx' \tens yy'), 
$$
for homogeneous element $x,x',y,y' \in \Uqgm$. 

Let $$\Delta_\n\cl \Uqgm \to \Uqgm \tens \Uqgm$$  be the algebra homomorphism defined by  
\begin{align} \label{Eq: Delta_n}
\Delta_\n(f_i) = f_i \tens 1+1 \tens f_i \quad \text{ for any $i \in I$.}    
\end{align} 
Note that
\begin{align} \label{eq: f_in pair}
\Delta_\n(f_i^{(n)}) = \sum_{k=1}^n q_i^{-k(n-k)} f_i^{(k)} \tens  f_i^{(n-k)} \quad \text{ for any $i\in I$ and $n \in \Zp$.}
\end{align}

The \emph{quantum unipotent coordinate ring} $\Aqn$ of $U_q(\g)$ is defined by 
$$
\Aqn \seteq \bigoplus_{\be \in \rl^-} \Aqn_\be, \quad \text{ where } \Aqn_\be \seteq \Hom_{\Q(q)}(\Uqgm,\Q(q)).
$$

Let  
$$
\bAng{ \ , \ } \cl \Aqn \times \Uqgm \to \Q(q)
$$
be the pairing.
We define the multiplication of $\Aqn$ as follows: for any $x,y \in \Aqn $,  
$$
\Ang{x y, a} \seteq    \bAng{x\tens y,\Delta_\n(a)} = \sum_{(a)} x (a_{(1)}) y( a_{(2)} ) \quad \text{ for any $a \in \Uqgm$,}
$$
where $\Delta_\n(a) = \sum_{(a)} a_{(1)} \otimes a_{(2)}$ is written in \emph{Sweedler notation}. We write $\wt(x) = \beta$ for any $x \in \Aqn_\beta$. 

For each $i \in I$, we denote by $\ang{i} \in \Aqn_{-\al_i}$ the dual element of $f_i$ with respect to $\bAng{ \ , \ }$; i.e., 
$$
\bAng{\ang{i},f_j} =\delta_{i,j} \quad \text{ for any $i,j\in I$.}
$$
Then the set $\{ \ang{i} \ | \ i\in I \}$ generates $\Aqn$. We set
$$
\Aqhn \seteq \Q(q^{1/2}) \tens_{\Q(q)} \Aqn. 
$$

\subsection{Bilinear form on $\Uqgm$} In this subsection, we review the symmetric bilinear form on $\Uqgm$.

Let us define a bilinear form $\kform{ \ , \ }$ on $\Uqgm$ by 
\begin{subequations} \label{eq: K form}
\begin{gather}
\kform{xy,z}  = \kform{x\tens y,\Delta_\n(z)} \quad \text{for any $x,y,z \in \Uqgm$}, \label{eq: K form1}\\
\kform{1,1}=1\qtq  \kform{f_i,f_i}  = 1 \quad \ \ \ \ \text{for any $i \in I$}, \label{eq: K form2}\\
\kform{x,y } =0 \quad \text{if  $x \in \Uqgm_\al$ and $y \in \Uqgm_\beta$ such that $\al \ne\beta$.}
\end{gather}
\end{subequations}
Then one can show that $\kform{ \ , \ }$  is unique, symmetric and non-degenerate by using the arguments in  \cite[Chapter 1]{LusztigBook}.

Let $e_i'$ (resp.\ $\es_i$) be the adjoint of the left (resp.\ right) multiplication of $f_i$ with respect to $\kform{ \ , \ }$; i.e.,
$$
\kform{ e_i'(x) , y } = \kform{ x,f_iy } \qtq \kform{ \es_i(x) , y } = \kform{ x,yf_i } \quad \text{ for any } x,y \in \Uqgm 
$$
(see \cite{K91,LusztigBook} and also~\cite[Section 2]{Kimura12}).  It is known that those adjoints are characterized by 
$$
[e_i,x] = \dfrac{\es_i(x)t_i - t_i^{-1}e'_i(x)}{q_i-q_i^{-1}} \quad\text{for $x \in \Uqgm$,}
$$
and satisfy 
$$
\Deltan(x) \equiv x \tens 1 + 1 \tens x + \sumI\left(  f_i \tens e'_i(x) + \es_i(x) \tens f_i  \right) \ \ {\rm mod} \; A
 \ \text{ for any }x \in\sumI \Uqgm f_i
$$
where $A = \sum_{\het(\al),\het(\be)>1} \Uqgm_\al \tens \Uqgm_\be$. For any $i,j \in I$ and $x,y \in \Uqgm$, we have
$$
e_i'\es_j =\es_j e'_i, \quad \es_i = * \circ e_i' \circ * \qtq \kform{x^*,y^*} = \kform{x,y}
$$
and 
\begin{subequations} \label{eq: e' es}
\begin{gather}
e'_i(xy)  = e'_i(x)y + q^{(\al_i,\wt(x))} x e'_i(y), \label{eq: e' 1}\\    
\es_i(xy) = x\es_i(y) + q^{(\al_i,\wt(y))} \es_i(x)y.   \label{eq: es 1} 
\end{gather}    
\end{subequations}
In particular, we have
\begin{align} \label{eq: e' f q-commute}
e_i'f_j = q^{ -(\al_i,\al_j) } f_je'_i +\delta_{i,j}.    
\end{align}

Note that
$$
\bigcap_{i\in I} \Ker \; e'_i = \bigcap_{i\in I} \Ker \; \es_i = \Q(q).
$$

\begin{remark}
The bilinear form $\kform{ \ , \ }$ on $\Uqgm$ given in \eqref{eq: K form} coincides with Kashiwara's bilinear form on $\Uqgm$	introduced in \cite[Proposition 3.4.4]{K91}.
The relation with Lusztig's bilinear form $(\ , \ )_L$ (\cite[Proposition 1.2.3]{LusztigBook}) is described in \cite[Section 2.2]{L04}. 
\end{remark}

\vskip 2em

\section{Quantum boson algebras} \label{Sec: QBA}
In this section, we investigate several properties of quantum boson algebras.

\begin{definition}
The \emph{quantum boson algebra} $\Bqg$ is the subalgebra of $\End_{\Q(q)}(\Uqgm)$ generated by the left multiplications of $f_i$ and $e'_i$
for all $i \in I$. 
\end{definition}

Let 
\eq
\sfE\cl \Uqgm\to\Bqg\label{Eq:sfE}
\eneq
be the $\Qq$-algebra homomorphism given by $f_i\mapsto e'_i$ for $i\in I$ and
let
\begin{align} \label{Eq: sfL}
\sfL \cl \Uqgm\to\Bqg
\end{align}
be the $\Qq$-algebra homomorphism such that $\sfL(x)$ is the left multiplication by $x\in\Uqgm$.
Then we have
$$
\kform{\sfE(x)y,z} = \kform{y,x^*z} \quad \text{ for any $x,y,z \in \Uqgm$.}
$$

We define a $\Q(q)$-linear map
\begin{align} \label{eq: Psi}
    \Psi\cl \Uqgm \tens \Uqgm \to \Bqg
\end{align}
defined by $x \tens y \mapsto \sfE(x)\sfL(y)$ for $x ,y \in \Uqgm$. It is obvious that $\Psi$ is surjective by~\eqref{eq: e' f q-commute}. 

\begin{proposition} \label{prop: bijection}
The $\Q(q)$-linear map $\Psi$ in~\eqref{eq: Psi} is bijective.      
\end{proposition}

\begin{proof}
For $  y  \in \Uqgm$,~\eqref{eq: e' es} says that
$$
\es_i \circ \sfL(y) = \sfL(y) \circ  \es_i + \sfL(\es_i(y)) Q_i 
$$
where $Q_i \in \End_{\Q(q)}(\Uqgm)$ is define by $Q_i(v) = q^{(\al_i,\wt(v))}v$ for a homogeneous element $v \in \Uqgm$. 

Since $e_j'\es_i = \es_i e'_j$, we have
$$
\es_i \circ \Psi(x\tens y) = \Psi(x\tens y) \circ \es_i + \Psi(x\tens \es_i y) Q_i 
$$
for elements $x,y \in \Uqgm$. Hence, for any $z \in \Uqgm \tens \Uqgm$, we obtain 
\begin{align} \label{eq: S_i}
\es_i \circ \Psi(z) =\Psi(z) \circ \es_i + \Psi(S_i(z))Q_i     
\end{align}
where $S_i= {\rm id}_{\Uqgm} \tens \es_i \in \End(\Uqgm \tens \Uqgm)$. 

Now we are ready to prove the injectivity of $\Psi$. For $m \in \Zp$, let us set
$$
K_m \seteq \sum_{\het(\be)\le m} \Uqgm \tens \Uqgm_{\be}. 
$$
We shall show that $K_m \cap \Ker(\Psi) =0$ by induction on $m \in \Zp$. The case $m=0$ is trivial. Now let us assume that $m>0$.
For any $i \in I$ and $z \in K_m \cap \Ker(\Psi)$,~\eqref{eq: S_i} implies that
$$
0 = \es_i \circ \Psi(z) = \Psi(S_i(z))Q_i.
$$
Since $S_i(z) \in K_{m-1}$, we have
$$S_i(z) \in K_{m-1} \cap \Ker(\Psi)=0.$$
Hence we have
\begin{align*}
z \in \sbcap_i \Ker(S_i) & = \Uqgm \tens \big(\sbcap_i \Ker(\es_i Q_i) \big) = \Uqgm \tens \big(\sbcap_i \Ker(\es_i) \big)     \\
& = \Uqgm \tens \Q(q) = K_0,
\end{align*}
which implies that $z=0$. 
\end{proof}

\begin{corollary} \label{cor: Bqg relation}
The $\Q(q)$-algebra $\Bqg$ is the $\Q(q)$-algebra generated by $\{ e'_i ,f_i \ |  \ i\in I\}$ subject to following relations:
\begin{equation} \label{eq: presentation}
\begin{aligned}
& e'_i f_j = q^{-(\al_i,\al_j)} f_j e'_i +\delta_{i,j}, \\
&\sum_{k=0}^{b_{ij}} (-1)^k  \begin{bmatrix} b_{ij} \\ k     \end{bmatrix}_i  f_i^k  f_j f_i^{b_{ij} - k} 
= \sum_{k=0}^{b_{ij}} (-1)^k  \begin{bmatrix} b_{ij} \\ k     \end{bmatrix}_i e_i^{\prime k }e'_j e_i^{\prime b_{ij} - k} =0 \ \  \text{ for }i \ne j,
\end{aligned}
\end{equation}
where $b_{ij}=1 - \lan h_i,\al_j \ran.$
\end{corollary}

\begin{proof}
Let $B$ be the $\Qq$-algebra generated by
$\{e'_i, f_i\mid i\in I \}$ with the same defining relations as above.
We then have a surjective algebra homomorphism
$$ 
a \cl B\epito\Bqg
$$ sending $e_i'$, $f_i$ to $\sfE(f_i)$, $\sfL(f_i)$  respectively by the $q$-Serre relations and \eqref{eq: e' f q-commute}.

On the other hand, let $E\cl \Uqgm\to B$ and $L \cl \Uqgm\to B$ be the homomorphisms  given by $f_i\mapsto e'_i$ and $f_i\mapsto f_i$ respectively.
 We then have the linear map 
$$
b \cl \Uqgm\tens\Uqgm\to B
$$  defined by
$b(x\tens y)=E(x)L(y)$.
Since
the composition $a \circ b$ coincides with $\Psi$ given in \eqref{eq: Psi}, 
$$
\xymatrix{
\Uqgm\tens\Uqgm  \ar@/_{-2.0pc}/[rrrr]^{\Psi} \ar@{->>}[rr]^{\qquad \quad b} &&  B \ar[rr]^{a \qquad } && \Bqg
}
$$
Proposition \ref{prop: bijection} says that $B\to\Bqg$ is bijective.
\end{proof}

\vskip 2em

\section{Upper global basis} \label{Sec: ugb}

Note that there exists a $\Q(q)$-algebra isomorphism 
\begin{align} \label{eq: iota}
\iota\cl \Uqgm \isoto \Aqn     
\end{align}
defined by $\iota(f_i) \seteq (1-q_i^2)^{-1}\ang{i}$ for any $i\in I$. We define a bilinear form $\aform{ \ , \ }$ on $\Aqn$ by 
\begin{align} \label{eq: aform}
\aform{f,g} \seteq \bAng{ f, \iota^{-1}(g)} \quad \text{ for any } f,g \in \Aqn.     
\end{align}

For any $i,j \in I$, we have
\begin{align} \label{eq: aform property}
\aform{\ang{i},\ang{j}} =\delta_{i,j} (1-q_i^2) \qtq \aform{\ang{ij},\ang{ij}} = \dfrac{(1-q_i^2)(1-q_j^2)}{1-q^{-2(\al_i,\al_j)}} \ \ \text{ if } (\al_i,\al_j)<0,    
\end{align}
where $\ang{ij} = \dfrac{\ang{i} \ang{j} - q^{-(\al_i, \al_j)   }\ang{j} \ang{i}   }{ 1-q^{-2(\al_i, \al_j)} }$ when $ (\al_i, \al_j)<0$. 
Hence we have
$$\text{$\bAng{\ang{ij},f_if_j}=1$ and $\bAng{\ang{ij},f_jf_i}=0$.}$$

For any $i \in I$, $\ep \in  \st{+,-} $  and $ \be  = \sumI n_i\al_i \in \rl$, we set
\begin{align} \label{eq: zeta}
\zeta_{i,\ep} = 1 -q_i^{\ep 2} \qtq \zeta^\be_\ep \seteq \prod_{i\in I} \zeta_{i,\ep}^{n_i} = \prod_{i\in I} (1-q_i^{\ep 2})^{n_i}.     
\end{align}
 We usually use $\zeta_i$  for $\zeta_{i, +}$  and $\zeta^\be$  for $\zeta_+^\be$.
Similarly we set
\eq \qq^\be \seteq \prod_{i\in I} q_i^{n_i}.
\label{eq:qq}
\eneq
Hence we have
$$    \ubzeta[\beta]  =(-1)^{\het(\beta)}\zeta^\beta\qq^{-2\beta}.$$

It is easy to see the following lemma.

\begin{lemma} \label{lem: zeta}
The bilinear form $\aform{\ ,\ }$ on $\Aqn$ has the following properties.
  \bnum
\item $\aform{\ ,\ }$ is a non-degenerate symmetric bilinear form
on $\Aqn$. 
  \item
    Let $u,v \in \Uqgm$  and set $\beta \seteq  \wt(u)\in\rl^-$. Then we have
    \begin{align*} 
\ang{\iota(u),v}=\aform{  \iota(u), \iota(v) }=\zeta^{\beta} \kform{u,v},
\end{align*}	
where $\iota\cl\Uqgm \buildrel \sim \over \longrightarrow \Aqn $ is given in  \eqref{eq: iota}. 
\ee
\end{lemma}

Let $\overline{\phantom{a}}$ be the $\Q$-algebra anti-automorphism on $\Uqgm$ defined by 
$$
\overline{q}=q^{-1} \qtq \overline{f_i} =f_i. 
$$
Define the map $\sfc \in \End(\Aqn)$ by
\begin{align} \label{eq: sfc}
 \bAng{\sfc(f),x} = \overline{\bAng{f,\overline{x}}} \quad \text{ for any $f\in \Aqn$ and $x \in \Uqgm$.}    
\end{align}
Then it satisfies 
\begin{align} \label{eq: property of sfc}
\sfc(q) =q^{-1}  \qtq \sfc(fg) = q^{(\wt(f),\wt(g))} \sfc(g)\sfc(f) \quad \text{ for } f,g \in \Aqn   
\end{align}
(see~\cite[Proposition 3.6]{Kimura12} for  example). 

We define 
\begin{align} \label{eq: Aqn integral form}
\Azn \seteq \{ f \in \Aqn \ | \ \bAng{f,U_\Zq^-(\g)} \in \Zq \}.    
\end{align}
Then $\Azn$ is a $\Zq$-subalgebra of $\Aqn$, since
$$\Delta_\n(\Uzgm)\subset\Uzgm\tens_{\Zq}\Uzgm.$$ 
Set 
$$\bbA \seteq \Z[q^{\pm1}, (1-q^n)^{-1} \ | \  n\in \Zsp \} \subset \Q(q).$$
Then we have the isomorphism
$$
\bbA \tens_{\Zq} U^-_\Zq(\g) \isoto  \bbA \tens_{\Zq}  A_\Zq(\n)
$$
induced by the isomorphism $\iota$ in~\eqref{eq: iota} (\cite[Theorem 5.7]{Tani}). Thus $\bbA \tens_{\Zq} A_\Zq(\n)$
is the $\bbA$-subalgebra of $\Aqn$ generated by $\st{ \ang{i}\mid i \in I }$.

\begin{lemma} \label{Lem: mem for An and ei'}
  Assume that $   a   \in \Um$ satisfies
  $\iota(a)\in\Azn$. Then, for any $n\in \Z_{\ge0}$, we have
	$$
	\iota \left( \zeta_i^{-n} e_i'^{(n)} (a) \right) \in \Azn.
	$$	
\end{lemma}
\begin{proof}
 
	Let $b \in\Uzgm$ and set $\beta=\wt(b)$.
	Since $ f_i^{(n)} b  \in \Uzgm$, it follows from Lemma~\ref{lem: zeta}
        that
	\begin{align*}
		\ang{ \iota ( \zeta_i^{-n} e_i'^{(n)} (a) , b  } 
		& =  \zeta^{ \beta } \zeta_i^{-n} (   e_i'^{(n)} (a)  ,  b) \\
		& =  \zeta^{\beta - n\al_i } (    a , f_i^{(n)}  b ) \\
		& =  \ang{\iota(a), f_i^{(n)} b} \in \Zq,
	\end{align*}	
	which gives the assertion.
\end{proof}

Let
$\{ \Gup(b) \ | \ b \in B(\infty) \}$ be the upper global basis of $\Azn$ (see \cite{K91, K93, K95} for its definition and properties). Set
$$
\LupA \seteq \sum_{b \in B(\infty)} \Z[q]\Gup(b) \subset \Azn.
$$

We regard $B(\infty)$ as a basis of $\LupA/q\LupA$ by 
\begin{align} \label{eq: local basis}
b \equiv \Gup(b) \qquad {\rm mod} \; q\LupA.    
\end{align}

We know that $\aform{\Gup(b),\Gup(b')}|_{q=0} =\delta_{b,b'}$ and
hence $B(\infty)$ is an orthonormal basis of
$\LupA/q\LupA$, which implies that
the lattice $\LupA$ is characterized by 
$$
\LupA = \{ x \in \Azn \ | \  \aform{x,x} \in \Z[[q]] \subset \Q(\hspace{-.4ex}(q)\hspace{-.4ex}) \}.
$$
We have
\begin{align} \label{eq: a-value}
\Aform{\LupA,\LupA} \subset \Z[q,(1-q^n)^{-1} \ | \ n\in \Zsp]    
\end{align}
(see \cite[Theorem 5.7]{Tani}). 

\vskip 2em 

\section{Bosonic Extensions} \label{Sec: BE}

In this section, we introduce and investigate a bosonic extension of \emph{arbitrary symmetrizable Kac-Moody type}. 
This is a  generalization of the result of \cite{OP24}
which treats bosonic extensions of finite type. 
We basically follow the notations in \cite{OP24}.  
Let $\sfC=(c_{i,j})_{i,j\in I}$ be a generalized Cartan matrix of symmetrizable Kac-Moody type.

\smallskip

We define $\hcalA$ to be the $\Qh$-algebra generated by the generators
$\{f_{i,p}\mid i\in I,p\in \Z\}$ with the defining relations:
\eq \label{Eq: def of hA}
&&\hspace{1ex}\left\{
\parbox{75ex}{
\bna
\item \label{it: def of hA (a)}
$\sum\limits_{k=0}^{1- \langle h_i,\al_j \rangle}(-1)^k\bnom{1-\Ang{h_i,\al_j}\\k}_i
f_{i,p}^kf_{j,p}f_{i,p}^{1-\langle h_i,\al_j \rangle-k}=0$\quad for any $i,j\in I$ such that $i\not=j$ and $p\in \Z$,

\vs{1ex}
\item \label{Eq: def of hA (b)}
$f_{i,m}f_{j,p}=q^{(-1)^{p-m+1}(\al_i,\al_j)}f_{j,p}f_{i,m}
+\delta_{(j,p),\,(i,m+1)}(1-q_i^2)$\quad if $m<p$,
\ee
}\right.
\eneq
where $ \delta_{(j,p),\,(i,m+1)} = \delta \left( (j,p) = (i,m+1) \right)$.

We set
\begin{align} \label{eq: wt of hcalA}
\al_{i,m} = (-1)^m \al_i \quad \text{ for $i \in I$ and $m \in \Z$}.     
\end{align}
Then~\eqref{Eq: def of hA (b)} reads as
$f_{i,m}f_{j,p}=q^{-(\al_{i,m},\al_{j,p})} f_{j,p} f_{i,m} +\delta_{(j,p),\,(i,m+1)}(1-q_i^2)$ if $m<p$. In particular,
$$
f_{i,m}f_{i,m+1} = q_i^2 f_{i,m+1}f_{i,m} + (1-q_i^2).  
$$
For any $n\in \Zp$, we set $f_{i,m}^{(n)} \seteq f_{i,m}^n/[n]_i!$. 

For $\be \in \rl$, let $\hcalA_\be$ be the $\Qh$-vector subspace of $\hcalA$ generated by
$f_{i_1,p_1}\cdots f_{i_r,p_r}$ such that $r \in \Zp$, $i_k \in I$, $p_k\in\Z$ and
$-\sum_{k=1}^r \al_{i_k,p_k}=\be$. Then we have the weight space decomposition
$$
\hcalA = \soplus_{\be \in \rl} \hcalA_\be.
$$
We say that an element $x \in \hcalA_\be$ is homogeneous of weight $\beta$
and set 
$$
\wt(x) = \be. 
$$ 
Note that $\wt(f_{i,p})=-\al_{i,p}=(-1)^{p+1}\al_i$.

 By the definition of $\hcalA$, we have  the following (anti-)automorphisms (see \cite[Section 3]{OP24}):
\bna
\item the anti-automorphism of $\Qh$-algebra 
${}^\star\cl \hcalA \longrightarrow \hcalA$ defined by 
\begin{align} \label{eq: star}
(f_{i,p})^\star = f_{i, -p},
\end{align}
\item the anti-automorphism of $\Q$-algebra $ \calD \cl \hcalA \longrightarrow \hcalA$ defined by 
$$
\calD (q^{\pm1/2}) = q^{\mp1/2}, \qquad \calD(f_{i,p}) = f_{i, p+1}, 
$$ 
\item  the anti-automorphism of $\Q$-algebra $\overline{\phantom{a}} \cl \hcalA \longrightarrow \hcalA$ defined by
$$
\overline{q^{\pm1/2}}=q^{\mp1/2}\qtq \overline{f_{i,p}}=f_{i,p},
$$ 
\item the automorphism of $\Qh$-algebra 
  $ \ocalD \seteq \overline{\phantom{a}} \circ \calD  = \calD \circ \overline{\phantom{a}}$ of $\hcalA$ given by\renewcommand{\thefootnote}{\ensuremath{\fnsymbol{footnote}}}
\footnote{Note that $\ocalD$ in our paper is denoted by $\calD$ in \cite{OP24}.}
\begin{align} \label{Eq: shift}
  \ocalD(f_{i,p}) = f_{i, p+1} \qquad \text{ for any $i\in I$ and $p\in \Z$.}
\end{align}
\ee

We define a map $c\cl \hcalA \to \hcalA$ by 
\begin{align} \label{eq: c map}
c(x)\seteq q^{N(\wt(x))} \overline{x} \quad\text{for any homogeneous element } x \in \hcalA,
\end{align}
where we set 
\begin{align} \label{eq: N}
N(\al) \seteq (\al,\al)/2 \quad\text{for }\al\in \rl.    
\end{align}
It is obvious that $N$ is a $\Z$-valued quadratic form on $\rl$ and satisfies 
\begin{align} \label{eq: N-properties}
N(\al+\be) = N(\al)+N(\be)+(\al,\be) \quad\text{for any } \al,\be\in \rl.    
\end{align}

By the definition, we have
\begin{align} \label{eq: c anti}
c(q)=q^{-1} \qtq c(xy) = q^{(\wt x,\wt y)} c(y)c( x )  \quad\text{for any homogeneous } x,y \in \hcalA.    
\end{align}

We set 
\begin{align} \label{eq: sigma}
\sigma(x) \seteq q^{-N(\wt\,x)/2} x \quad\text{for any homogeneous element } x \in \hcalA.    
\end{align}

Note that $N(\wt\,x)/2\in\frac{1}{2}\Z$.  It is easy to see that
$$
\sigma \circ c = \overline{\phantom{a}} \circ \sigma ,
$$
and $\sigma$ sends $c$-invariant elements to bar-invariant elements. Namely, we have
\eq \label{Eq: bar inv c inv}
&&\hs{-5ex}
\ba{l}
\bullet\ \text{if $x$ is bar-invariant, then $ \sigma^{-1}(x) =  q^{(\wt(x),\wt(x))/4}x$ is $c$-invariant,}\\
\bullet\ \text{if $x$ is $c$-invariant, then $\sigma(x) = q^{-(\wt(x),\wt(x))/4}x$ is bar-invariant.}
\ea
\eneq

Note that
$$
\sigma(xy) = q^{-(\wt( x),\wt( y))/2} \sigma(x)\sigma(y) \qtq \sigma(x^{n}) = q^{-n(n-1)(\wt (x),\wt (x))/4} \sigma(x)^n.
$$

\begin{definition}
For $ - \infty \le a \le b \le \infty$, let $\hcalA[a,b]$ be the $\Q(q^{1/2})$-subalgebra
of $\hcalA$ generated by $\{ f_{i,k} \ | \  i \in I, a\le k \le b\}$.
We simply write 
$$
\hcalA[m] \seteq \hcalA[m,m], \quad 
\hcalA_{\ge m} \seteq \hcalA[m,\infty], \quad
\hcalA_{\le m} \seteq \hcalA[-\infty,m]. 
$$
Similarly, we set $\hcalA_{>m} \seteq \hcalA_{\ge m+1}$
and $\hcalA_{<m} \seteq \hcalA_{\le m-1}$. 
\end{definition}

We define the weight spaces $\hcalA[a,b]_\be$ by $\hcalA[a,b] \cap \hcalA_\be$. Note that $\hcalA[k]_\be=0$ for $\be \not\in (-1)^{k+1}\rl^+$. 

In the rest of this section, we shall show that there is a realization of $\hcalA$.

For $m\in\Z$, let us denote by
\eq
\bfL_m\cl\Uqhm\to\hA[m]
\eneq
the $\Qh$-algebra homomorphism given by
$\bfL_m(f_i)=f_{i,m}$.

Let $\sfE_\bzeta\cl \Uqgm \to \Bqg $ be the $\Qq$-algebra homomorphism defined
by $\sfE_\bzeta(f_i) = \bzetai e'_i$ for $i \in I$.
Let $\Phi_{m}$ be the $\Qq$-algebra homomorphism 
\eq \label{eq: Phi bzeta}
&&\Phi_{m}\cl    \Bqg   \to \hcalA \quad
\qt{given by $e'_i \longmapsto   \ubzetai{i}{-1}   f_{i,m+1}$, 
\quad $f_i \longmapsto f_{i,m}$.}
\eneq
  It is well defined by Corollary~\ref{cor: Bqg relation}.

\begin{lemma} \label{lem: Ezeta}
For $i \in I$ and an homogeneous element $y \in \Uqgm$, we have
$$
\sfL(f_i) \sfE_\bzeta(y) = q^{-(\al_i,\wt(y))} \left( -\bzetai\sfE_\bzeta(\es_i(y)) + \sfE_\bzeta(y)\sfL(f_i) \right). 
$$
\end{lemma}

\begin{proof}
It follows from~\eqref{eq: e' 1} that
\begin{align} \label{eq: sfE sfLu}
\sfE(f_i) \sfL(u) = q^{(\al_i, \wt(u)) } \sfL(u)\sfE(f_i) + \sfL(e'_i(u))    
\end{align}
for any homogeneous $u \in \Uqgm$. Let $\fraks$ be the anti-automorphism of $\Bqg$ defined by $f_i \mapsto e'_i$ and $e'_i \mapsto f_i$ for $i \in I$. Note that
$\fraks (\sfE(u)) = \sfL(u^*)$ and $\fraks (\sfL(u)) = \sfE(u^*)$. Applying $\fraks$ to~\eqref{eq: sfE sfLu} when $u = y^{ *}$, we have
$$
\sfE(y) \sfL(f_i) = q^{(\al_i,\wt(y))} \sfL(f_i)\sfE(y) + \sfE(\es_i(y)),    
$$
which says 
$$
\sfL(f_i)\sfE(y)  = q^{-(\al_i,\wt(y))} ( \sfE(y)\sfL(f_i) - \sfE(\es_i(y)).    
$$
This gives the assertion by rewriting in terms of $\sfE_\bzeta$
 since $\sfE_\bzeta(y)= \ubzeta[-\wt(y)]  \sfE (y)$. 
\end{proof}

 Let $\Uqhm \seteq  \Q(q^{1/2}) \tens_{\Q(q)} \Uqgm $,
and for any $k\in\Z$, let $\Uqhm_k$ be the copy of $\Uqhm$.
For $a,b\in\Z$ such that $a\le b$, set
$$\scrU_\Z[a,b]=\Uqhm_b\tens\Uqhm_{b-1}\tens\cdots\tens\Uqhm_{a+1}\tens\Uqhm_a,$$
where  the tensoring is taken over $\Q(q^{1/2})$.
For $a\le a'\le b'\le b$, we embed
$\scrU_\Z[a',b']$  into $\scrU_\Z[a,b]$ by $u_{b'}\tens\cdots\tens u_{a'}
\mapsto u_b\tens\cdots\tens u_a$ with
$u_k=1$ for $k\in[a,a'-1]\cup[b'+1,b]$.
Then we set
$$\UZ=\sbcup_{a\le b}\UZ[a,b].$$

For  $ \bfu =  \displaystyle\tens_{k\in \Z}^{\to} u_k$ in $\UZ$ (with $u_k\in\Uqhm_k$
with $u_k=1$ for all but finitely many $k$'s),  we set 
$$
\bfu_{>m} \seteq \tens_{k \in \Z_{>m}}^{\to} u_k, \qquad 
\bfu_{<m} \seteq \tens_{k \in \Z_{<m}}^{\to} u_k,
$$
and define $\bfu_{\ge m}$ and $\bfu_{\le m}$ in a similar way. 
For $a<b$, we  have 
$$
\bfu = \bfu_{>b} \tens u_b \tens  u_{b-1} \tens \cdots \tens u_a \tens \bfu_{<a}.
$$
and define 
$$
\widehat{\wt}(\bfu) \seteq \sum_{k \in \Z} (-1)^k \wt(u_k) \quad\text{ for } \bfu \in \scrU_\Z.
$$
We also define $\hwt(\bfu_{>m})$, $\hwt(\bfu_{<m})$, etc. in a similar way.

Now let us define the action of $f_{i,m}$ $((i,m) \in I \times \Z)$ on $\scrU_\Z$ as follows: 
\begin{align} \label{eq: action f_im}
f_{i,m}(\bfu) \seteq q^{(\al_{i,m},\hwt(x\tens y))} (-\bzetai x \tens \es_i(y) \tens z \tens w  + x \tens y \tens  f_i z \tens w), 
\end{align}
where 
\begin{align} \label{eq: xyzw}
x = \bfu_{>m+1}, \quad  y = u_{m+1}, \quad z = u_m \qtq w = \bfu_{<m}. 
\end{align}

\begin{theorem} \label{Thm: repn of hA}
The action~\eqref{eq: action f_im} gives an $\hcalA$-module structure on $\scrU_\Z$.     
\end{theorem}

\begin{proof}
Let $i\in I$, $m \in \Z$ and $\bfu = \tens_{k \in \Z} u_k \in \scrU_\Z$. 

\smallskip

\noindent
{\rm (i)} We first deal with the quantum Serre  relation~\eqref{Eq: def of hA}~\eqref{it: def of hA (a)}.

Let
\begin{align}
\Psi_{\bzeta} \cl \Uqgm \tens \Uqgm \isoto \Bqg     
\end{align}
be the $\Q(q)$-linear map defined by $\Psi_\bzeta(x \tens y) = \sfE_\bzeta(x)\sfL(y)$.
Set $x,y,z$ and $w$ as in~\eqref{eq: xyzw}. By Lemma~\ref{lem: Ezeta}, the action~\eqref{eq: action f_im} can be written as 
$$
f_{i,m}(\bfu) = q^{(\al_{i,m},\hwt(x))} x \tens \left( \Psi_{\bzeta}^{-1}\big(\sfL(f_i)\Psi_{\bzeta}(y \tens z) \big) \right) \tens w. 
$$
Since the quantum Serre relations hold for $\sfL(f_i)\in\Bqg$, the same relations hold for $f_{i,m}\in\End(\scrU_\Z)$. 

\smallskip

\noindent
{\rm (ii)} We now consider the relation~\eqref{Eq: def of hA}~\eqref{Eq: def of hA (b)}. Let $m,p \in \Z$ with $p>m$. 

{\rm (a)} Suppose that $p>m+1$. Since
$$
\hwt((f_{j,p}\bfu)_{ \ge m+1}) =\hwt(\bfu_{ \ge  m+1}) -\al_{j,p} \qtq \hwt((f_{j,p}\bfu)_{\ge  p+1}) =\hwt(\bfu_{\ge p+1}) 
$$
by~\eqref{eq: action f_im}, we have 
$$
f_{i,m}f_{j,p}(\bfu) = q^{(\al_{i,m},-\al_{j,p})}f_{j,p}f_{i,m}(\bfu) = q^{(-1)^{p-m+1}(\al_i,\al_j)}f_{j,p}f_{i,m}(\bfu). 
$$

{\rm (b)} Suppose that $ p=m+1$. 
We set 
\begin{align*}
&a = \bfu_{> m+2}, \quad  b = u_{m+2}, \quad \  c = u_{m+1}, \quad  \ \  d = u_{m}, \quad \  e = \bfu_{< m}, \\
&\al = \hwt(a), \quad  \beta = \hwt(b) , \quad \gamma = \hwt(c), \quad \delta = \hwt(d), \quad \epsilon  = \hwt(e), 
\end{align*}	
and simply write $\bfu=abcde$ instead of $a \tens b \tens c \tens d \tens e$. Let
$$
C = (\al_{i,m},\al+\be+\ga) + (\al_{j,m+1},\al+\be).
$$

We then have
\begin{align*}
& q^{-C-(\al_i,\al_j)}  f_{i,m}f_{j,p} (abcde) \\ 
& = \bzetai\bzetai[j] a (\es_j(b))(\es_i(c))de -\bzetai[j]  a(\es_j(b)) c (f_i d) e - \bzetai ab (\es_if_j(c)) de + ab(f_jc)(f_id)e,  
\end{align*}
and
\begin{align*}
& q^{-C}  f_{j,p}f_{i,m} (abcde) \\ 
& = \bzetai\bzetai[j] a (\es_j(b))(\es_i(c))de -\bzetai[j]  a(\es_j(b)) c (f_i d) e - \bzetai ab (f_j\es_i(c)) de + ab(f_jc)(f_id)e.  
\end{align*}
Thus we have
$$
(f_{i,m}f_{j,m+1}-q^{(\al_i,\al_j)}f_{j,m+1}f_{i,m}) (\bfu) = -q^{C+(\al_i,\al_j)}
\bzetai ab( (\es_if_j - f_j\es_i)c )de. 
$$

In the case when $i \ne j$, we have $\es_i f_j - f_j \es_i =0$. Thus it implies that
$$f_{i,m}f_{j,m+1}-q^{(\al_i,\al_j)}f_{j,m+1}f_{i,m} =0.$$ 

Suppose that $i =j$. Since $(\es_i f_j - f_j \es_i)(c) =q^{(\al_i,\wt(c))}c$ and 
$C =- (\al_i,\wt(c))$,  we have
$$
(f_{i,m}f_{j,m+1}-q^{(\al_i,\al_j)}f_{j,m+1}f_{i,m}) (\bfu) = -q_i^2\bzetai \bfu = (1-q_i^2)\bfu,
$$
as we desired. 
\end{proof}

\begin{corollary} \label{cor: main1} \hfill
\bnum
\item \label{it: alg homo}
  For each $m\in \Z$, the $\Qh$-algebra homomorphism $\bfL_m\cl \Uqhm \to
  \hcalA[ m]$ given by $f_i \mapsto f_{i,m }$ is an isomorphism.
\item For any $a,b \in \Z$ with $a \le b$, the $\Qh$-linear map 
$$
\hcalA[b] \tens_\Qh \hcalA[b-1] \tens_\Qh \cdots \tens_\Qh \hcalA[a+1] \tens_\Qh \hcalA[a] \to \hcalA[a,b]
$$
defined by $x_b \tens x_{b-1} \tens \cdots \tens x_{a+1} \tens x_a \longmapsto x_bx_{b-1}\cdots x_{a+1}x_a$ is an isomorphism. 
\ee
\end{corollary}

\begin{proof}
Let $\mathsf{1} \seteq \tens_{k \in \Z} 1 \in \scrU_\Z$ and define a $\Qh$-linear map 
$$  F \cl \hcalA \longrightarrow \scrU_\Z  \quad \text{ by } x \longmapsto x(\mathsf{1})$$
induced  by  the module structure~\eqref{eq: action f_im}, and define  
$$
G \cl \scrU_\Z \longrightarrow \hcalA \quad \text{ by } \bfu = \tens_{k \in \Z}^{\to} u_k \longmapsto \oprod_{k \in \Z} \bfL_k (u_k). 
$$
It is easy to see that $F \circ G = {\rm id}_{\scrU_\Z}$
 and $G\circ F={\rm id}_{\hA}$,  which gives the desired properties.
\end{proof}

Note that Corollary \ref{cor: main1} tells us that  $\hcalA$-module  $\scrU_\Z$ is faithful.

\begin{remark}
When $\g$ is of finite type, the action of $f_{i,k}$ on the regular representation is described in \cite[Theorem 6.4]{OP24}.
\end{remark}	

{}From Corollary $\ref{cor: main1}$, any element $x$ in $ \hcalA[a,b]$ can be expressed as 
$$
x = \sum_{t} x_{b,t} x_{b-1,t} \cdots x_{a,t},
$$  
where $x_{k,t}$ is a homogeneous element in $\hcalA[k]$ and $t$ runs over a finite set. 
An element $x \in \hA[a,b]$ is said to be \emph{strongly homogeneous} 
if $x$ can be expressed as $x = x_{b} x_{b-1} \cdots x_{a}$ for some homogeneous elements $x_k\in \hcalA[k]$. 
For a strongly homogeneous element $x = x_{b} x_{b-1} \cdots x_{a}$, we define
(see \eqref{eq: het dVert}) 
\begin{align} \label{Eq: Ht}
\het(x) \seteq \sum_{t=a}^b \het\bl\wt(x_{k})\br \in \Zp.
\end{align}

\vskip 2em 

\section{Bilinear forms on $\hcalA$} \label{Sec: hAform}

For homogeneous elements $x,y \in \hcalA$, we set
\begin{align*}
[x,y]_q \seteq xy - q^{-(\wt\, x,\wt\, y)}yx.     
\end{align*}

It is easy to see that
\begin{subequations} \label{eq: []q}
\begin{gather}
    [x,yz]_q = [x,y]_qz +q^{-(\wt\, x,\wt\, y)}y[x,z]_q, \label{eq: []q1}\\
    [xy,z]_q = x[y,z]_q +q^{-(\wt\,   y ,\wt\, z)}[x,z]_qy, \label{eq: []q2}
\end{gather}    
\end{subequations}
for any homogeneous elements $x,y,z \in \hcalA$. For any $i \in I$ and $m\in \Z$, let $\rmE_{i,m}$ and $\Es_{i,m}$ to be the endomorphisms 
of $\hcalA$ defined by
\footnote{
The operators $\rmE_{i,m}$ and $\Es_{i,m}$ are different from $E_{i,m}$ and $E_{i,m}^\star$ given in \cite[Section 6]{OP24}.}
\begin{align} \label{eq: Ei Esi}
\rmE_{i,m}(x) \seteq [x,f_{i,m+1}]_q \qtq  \Es_{i,m}(x) \seteq [f_{i,m-1},x]_q     
\end{align}
for any homogeneous element $x \in \hcalA$. For any $n\in \Zp$, we set
$$
\rmE_{i,m}^{(n)} \seteq \dfrac{1}{[n]_i!} \rmE_{i,m}^n, \qtq \Esn{(n)}_{i,m} \seteq \dfrac{1}{[n]_i!} \Esn{n}_{i,m}.
$$
Hence we have
\begin{subequations} \label{eq: E and Star}
\begin{gather}
\Es_{i,m} = \star \circ \rmE_{i,-m} \circ \star, \\
\rmE_{i,m}(f_{j,m}) = \Es_{i,m}(f_{j,m}) =\delta_{i,j}(1-q_i^2),
\end{gather}
\end{subequations}
where $\star$ is the anti-automorphism in~\eqref{eq: star}, and 
\begin{subequations} \label{eq: x f_im fim x}
\begin{gather}
x f_{i,m+1} =  \rmE_{i,m}(x) + q^{-(\al_{i,m},\wt\, x)} f_{i,m+1}x, \\
f_{i,m-1}x  =  \Es_{i,m}(x) + q^{-(\al_{i,m},\wt\, x)} x f_{i,m-1}.
\end{gather}
\end{subequations}
Note that $\rmE_{i,m}$ and $\Es_{i,m}$ have weight $\al_{i,m}$. 

For each $m \in \Z$, we define a $\Qh$-algebra isomorphisms 
\begin{align} \label{eq: vph}
	\varphi_{m}\cl: \Aqhn  \longrightarrow \hcalA[m] \qquad \text{ by } \varphi_m(\ang{i}) = q_i^{1/2} f_{i,m} . 
\end{align}
Then, we can easily see that
\eq
\vphi_m\bl\iota(a)\br=\qq^{-\wt(a)/2}\zeta^{\wt(a)}\bfL_m(a)\qt{for any
  $a\in\Um$.}\label{eq:phiiota}
\eneq
where $\bfL_m$ is given in Corollary~\ref{cor: main1} and $\qq$ is given in \eqref{eq:qq}.

\begin{lemma} \label{lem: E Es} Let $i \in I$ and $m \in \Z$.
\bnum
\item \label{it: E Es}
For any homogeneous $x,y \in \hcalA$,
\begin{align*}
\rmE_{i,m}(xy) &= x \rmE_{i,m}(y) + q^{-(\al_{i,m},\wt\, y)} \rmE_{i,m}(x)y, \\    
\Es_{i,m}(xy) &=  \Es_{i,m}(x)y + q^{-(\al_{i,m},\wt\, x)} x\Es_{i,m}(y).    
\end{align*}
\item \label{it: E comm}
For any homogeneous elements $x \in \hcalA_{>m}$ and $y \in \hcalA_{<m}$, we have
\begin{align*}
f_{i,m-1} x = q^{-(\al_{i,m},\wt\, x)} x f_{i,m-1} \qtq     yf_{i,m+1}  = q^{-(\al_{i,m},\wt\, y)} f_{i,m+1}y.
\end{align*}
Namely $$\rmE_{i,m}(\hA_{<m})=0\qtq\Es_{i,m}(\hA_{>m})=0.$$
\item \label{it: E and bfL}
For a  homogeneous $u \in \Uqhm$, we have 
\begin{align*}
\rmE_{i,m}(\bfL_{m}(u) ) & = q_i^{-\ang{h_i,\wt(u)}-2}\zeta_i \bfL_{m}(e'_iu),    \\
  \rmE_{i,m}(\vphi_m \circ\iota(u) ) & = q_i^{1/2-\ang{h_i,\wt(u)}-2}\vphi_m \circ \iota(e'_iu),\\
  \Es_{i,m}(\bfL_{m}(u) ) & = q_i^{-\ang{h_i,\wt(u)}-2}\zeta_i \bfL_{m}(\es_iu),\\
  \Es_{i,m}(\vphi_m  \circ \iota(u)) & =
 q_i^{1/2-\ang{h_i,\wt(u)}-2}\vphi_m  \circ \iota(\es_iu).
\end{align*}

In particular, for any $x \in \hcalA[m]$, we have $\rmE_{i,m}(x),\Es_{i,m}(x) \in \hcalA[m]$.   

\item \label{it: En}
For $f \in \Aqhn$ with $\be \seteq \wt(f)$, we have
$$
\rmE_{i,m}^{(n)}(\varphi_m \circ\iota(f))
= q_i^{n/2-n\ang{h_i,\be}-n(n+1)} \varphi_m \iota(e_i^{\prime(n)}f) .
$$
\ee
\end{lemma}

\begin{proof}
  \eqref{it: E Es} follows from \eqref{eq: []q}. 
  \eqref{it: E comm} follows from the definition~\eqref{Eq: def of hA} of $\hcalA$. Let us focus on~\eqref{it: E and bfL}. 
We  have the following
commutative diagram:
$$
\xymatrix{
\Uqhm\akete[-1.3ex]  \ar[rr]^{\bfL_m} \ar@{>->}[d]_{\sfL} && \hcalA[m]\akete[-1.3ex] \ar@{>->}[d]  \\
\Bqg  \ar[rr]^{\Phi_{m}} && \hA[m, m+1], 	
}
$$
where $\sfL$ and $\Phi_m$ are given in~\eqref{Eq: sfL} and~\eqref{eq: Phi bzeta} respectively, and the right vertical map is the natural embedding.

Let $u \in \Uqhm$ and set $x \seteq \bfL_{m}(u)$. Then we have
\begin{align*}
\rmE_{i,m}(x) & = x f_{i,m+1} - q^{-(\al_{i,m},\wt\, x)}f_{i,m+1} x \\
& = (1-q_i^{-2}) \Phi_{m}( \sfL(u) e'_i - q^{-(\al_i,\wt\,u)} e'_i \sfL(u) ) \\
& \overset{*}{=} (1-q_i^{-2}) \Phi_{m}\bl -q^{-(\al_i,\wt\,u)} \sfL(e'_iu)\br \\
& =  q_i^{-\ang{h_i,\wt(u)}-2} \zeta_i\bfL_{m}(e'_iu),
\end{align*}
where $\zeta_{i}$ is given in~\eqref{eq: zeta}, and 
the equality $\overset{*}{=}$ follows from~\eqref{eq: e' 1}.
Hence we have obtained the first equality.
The second equality follows from the first and \eqref{eq:phiiota}. 

The equalities for $\Es_{i,m}$ follow from~\eqref{eq: E and Star}. 

\snoi 
\eqref{it: En}\  By\eqref{it: E and bfL} and the standard induction argument on $n$, we have
\eqn
\rmE_{i,m}^{(n)}(\varphi_m\iota(f)) && = [n]_i^{-1} \rmE_{i,m}(q_i^A  \varphi_m(\iota(e_i^{\prime (n-1)} f )   ) ) \allowdisplaybreaks\\
&& = q_i^{A+B}  \varphi_{m}(\iota(e_i^{\prime (n)} \fc )   ) ) \allowdisplaybreaks
\eneqn
Here $A = (n-1)/2 - (n-1)\Ang{h_i,\be}-n(n-1)$ and
$B=1/2-\ang{h_i,\beta+(n-1)\al_i}-2$.
Then $A+B=n/2-n\ang{h_i,\beta}-n(n+1)$. 
\end{proof}

{}From Corollary~\ref{cor: main1},  we have the decomposition
\begin{align} \label{eq: hcalA decomposition}
\hcalA = \soplus_{(\be_k)_{k \in \Z} \in \rl^{\oplus \Z}} \oprod_{k \in \Z} \hcalA[k]_{\be_k}.     
\end{align}
Define 
\begin{align} \label{eq: def of Mn}
\Mn\colon \hcalA\longepito \Qh
\end{align}
to be the natural projection
$ \hcalA\longepito \oprod_{k\in\Z}\hcalA[k]_{0}\simeq\Qh$
 arising from the decomposition \eqref{eq: hcalA decomposition}.

\begin{lemma}\label{Lem: Mfx = Mxf}
For any $x\in\hcalA$, $i\in I$ and $m\in\Z$,
we have
\begin{align} \label{eq:Mfx}
\Mn(f_{i,m-1}x)=\Mn(xf_{i,m+1}).    
\end{align} 
\end{lemma}

\begin{proof}
We may assume that $x = uvw$ where $u \in \hcalA_{>m}$, $v \in \hcalA[m]$
and $w \in \hcalA_{<m}$ are homogeneous elements.  If $\wt(x) + \wt(f_{i,m\pm1}) \ne 0$, then we have
$$
\Mn(f_{i,m-1}x)=\Mn(xf_{i,m+1}) =0
$$
by~\eqref{eq: def of Mn}. Thus we may assume that $\wt(u)+\wt(v)+\wt(w)=\wt(f_{i,m})$.  

Lemma~\ref{lem: E Es} and~\eqref{eq: x f_im fim x} say that
\begin{align*}
f_{i,m-1} x = q^{-(\al_{i,m},\wt\, u)} u \left( \Es_{i,m}(v)w + q^{-(\al_{i,m},\wt\, v)} v f_{i,m-1}w \right).    
\end{align*}
Since $\Mn(uv f_{i,m-1}w) =0$ by the definition~\eqref{eq: def of Mn}, we obtain
$$
\Mn(f_{i,m-1}x) = q^{-(\al_{i,m},\wt\,u)} \Mn( u \Es_{i,m}(v)w ).
$$

On the other hand, Lemma~\ref{lem: E Es} and~\eqref{eq: x f_im fim x} also say that
\begin{align*}
 x f_{i,m+1} = q^{-(\al_{i,m},\wt\,w)} \left( u \rmE_{i,m}(v)w + q^{-(\al_{i,m},\wt\, v)} u f_{i,m+1}vw \right).   
\end{align*}
Thus we have
$$
\Mn(xf_{i,m+1}) = \Mn(q^{-(\al_{i,m},\wt\,w)} u \rmE_{i,m}(v)w ),
$$
as in the case of $\Mn(f_{i,m-1}x)$. Hence, unless $\wt(u)=\wt(w)=0$ and $\wt(v)=\wt(f_{i,m})$, we have
$$
\Mn(f_{i,m-1}x)=\Mn(xf_{i,m+1}) =0.
$$

In the case when $\wt(u)=\wt(w)=0$ and $\wt(v)=\wt(f_{i,m})$, the element $v$ is contained in $\Qh f_{i,m}$ and  hence $\rmE_{i,m}(v)=\Es_{i,m}(v)$, 
which implies the assertion. 
\end{proof}

We define a bilinear form on $\hcalA$ as follows:
\begin{align} \label{eq: hA form}
\hAform{x,y} \seteq \Mn(x \ocalD(y)) \in \Qh \quad \text{ for any } x,y \in \hcalA,   
\end{align}
where $\ocalD$ is the automorphism of $\hcalA$ given in~\eqref{Eq: shift}.

\begin{lemma} \label{lem: hAform}
Let $i \in I$ and $m \in \Z$.
\bnum
\item \label{it: hAform 1}
For any $x\in\hcalA$, we have $\Mn\bl x^\star\br=\Mn\bl\ocalD(x)\br=\Mn(x)$.
\item \label{it: hAform 2}  If $x$ and $y$ are homogeneous elements such that $\wt(x)\not=\wt(y)$, then $\hAform{x,y}=0$.
\item \label{it: hAform 3}
For any $x,y \in \hA$, we have 
$$
\hAform{x,y} = \hAform{\ocalD(x), \ocalD(y)} = \hAform{ y^\star,x^\star}.
$$
\item \label{it: hAform 4} For any $x\in \hcalA_{\ge m}$ and $y\in \hcalA_{\le m}$, we have
  $$\Mn(xy)=\Mn(x)\Mn(y).$$
\item \label{it: hAform 5} For any $x,z \in \hcalA_{\ge m}$ and $y,w \in \hcalA_{<m}$, we have 
$$
\hAform{xy, zw} = q^{( \wt(y), \wt(z) )} \hAform{x,z} \hAform{y,w}.
$$
In particular, we have 
$$
\hAform{xy, w}=0 \qt{ when $\wt(x) \ne 0$.}
$$
\item \label{it: hAform 6} For any  $x,y\in \hA$, we have 
$$
\hAform{f_{i,m}x, y } =  \hAform{x, y f_{i,m+1} } \qtq  \hAform{xf_{i,m}, y } =  \hAform{x, f_{i,m-1} y }.
$$
\item \label{lem: properties of ( , ) (iv)}
For any $x,y\in\hcalA_{\le m}$ and $u, v \in\hcalA_{\ge m}$, we have
$$
\hAform{f_{i,m}x,y}=\hAform{ x,\rmE_{i,m}(y)} \qtq \hAform{u , vf_{i,m} }=\hAform{\Es_{i,m}(u) ,v}.
$$
\ee
\end{lemma}

\begin{proof}
\eqref{it: hAform 1}, \eqref{it: hAform 2} and \eqref{it: hAform 3} follow  from the definition \eqref{eq: def of Mn} of $\Mn$.
Let us only show $\hAform{x,y}=\hAform{y^\star,x^\star}$.
We have
\eqn
\hAform{x,y}&&=\Mn(x\ocalD y)=
\Mn\bl(x\ocalD y)^\star\br
=\Mn\bl(\ocalD y)^\star x^\star\br\\
&&=\Mn\bl\ocalD^{-1}(y^\star)x^\star\br
=\Mn\Bigl(\ocalD\bl\ocalD^{-1}(y^\star)x^\star\br\Bigr)
=\Mn\bl y^\star\ocalD(x^\star)\br
=\hAform{y^\star, x^\star}.
\eneqn

\snoi
\eqref{it: hAform 4}
We may assume that $x=x_{>m}x_m$ and $y=y_my_{<m}$ with
$x_{>m}\in\hcalA_{>m}$, $x_m\in\hcalA[m]$, $y_m\in\hcalA[m]$ and $y_{<m}\in\hcalA_{<m}$.
Then we have
$xy=x_{>m}(x_my_m)y_{<m}$ and hence we have
$$
\Mn(xy)=\Mn(x_{>m}(x_my_m)y_{<m})=
\Mn(x_{>m})\Mn(x_my_m)\Mn(y_{<m}).$$
Then the result follows from
$\Mn(x)=\Mn(x_{>m})\Mn(x_m)$,
$\Mn(y)=\Mn(y_m)\Mn(y_{<m})$ and $\Mn(x_my_m)=\Mn(x_m)\Mn(y_m)$.

\snoi
\eqref{it: hAform 5} It follows from 
\begin{align*}
\hAform{xy, zw} &= \Mn(xy\ocalD(zw)) =  	\Mn(xy\ocalD(z)\ocalD(w)) = 
 q^{-(\wt(y), \wt(\ocalD(z)))} \Mn(x\ocalD(z) y\ocalD(w)) \\
 &=  q^{(\wt(y), \wt( z))} \hAform{x,z} \hAform{y,w},
\end{align*}	
where the third follows from the second condition of \eqref{Eq: def of hA}, and the last identity follows from~\eqref{it: hAform 4},
$x\ocalD(z)\in\hcalA_{\ge m}$ and $y\ocalD(w)\in\hcalA_{\le m}$.

\mnoi
\eqref{it: hAform 6} By Lemma  \ref{Lem: Mfx = Mxf}, we have 
\begin{align*}
\hAform{f_{i,m}x, y } = \Mn(f_{i,m} x \ocalD(y) ) = \Mn(x \ocalD(y) f_{i,m+2}  ) =	 \Mn( x \ocalD(y f_{i,m+1}  ) ) =  \hAform{x, y f_{i,m+1}},
\end{align*}	
and 
\begin{align*}
\hAform{xf_{i,m}, y } = \Mn( xf_{i,m} \ocalD(y) ) = \Mn(x \ocalD(f_{i,m-1} y)  )  =  \hAform{x, f_{i,m-1} y }.
\end{align*}	

\mnoi
\eqref{lem: properties of ( , ) (iv)} It follows from~\eqref{it: hAform 6}, \eqref{eq: x f_im fim x} and~\eqref{it: hAform 3} that
\eqn
\hAform{f_{i,m}x,y}=\hAform{x,yf_{i,m+1}}
=\hAform{ x,\rmE_{i,m}(y)+q^{-(\al_{i,m},\wt y)}f_{i,m+1}y}
=\hAform{ x,\rmE_{i,m}(y)}.
\eneqn
Applying ${}^\star$ to $\hAform{f_{i,m}x,y}=\hAform{ x,\rmE_{i,m}(y)} $,
we obtain by \eqref{eq: E and Star} and (iii)
\eqn
\hAform{u, vf_{i,m}}
&&=\hAform{ (vf_{i,m})^\star,u^\star}
=\hAform{ f_{i,-m}v^\star,u^\star}\\
&&=\hAform{ v^\star,\rmE_{i,-m}(u^\star)}
=\hAform{ v^\star,(\Es_{i,m}u)^\star}
=\hAform{ \Es_{i,m}u, v}. \qquad \qedhere
\eneqn
\end{proof}

Now we present the main theorem of this section.

\begin{theorem} \label{thm: hAform}
Let $\hAform{ \ , \ }$ be the bilinear form on $\hcalA$ given in~\eqref{eq: hA form}.
\bnum
\item \label{it: hAform sym nond}
The bilinear form $\hAform{ \ , \ }$ is symmetric and non-degenerate.
\item \label{it: hAform decomposition}
For $x \in \oprod_{k \in [a,b]} x_k$ and $y \in \oprod_{k \in [a,b]} y_k$ with $x_k,y_k \in \hcalA[k]$, we have
$$
\hAform{x,y} = q^{\sum_{s<t} (\wt(x_s),\wt(x_t))} \prod_{k\in [a,b]} \hAform{x_k,y_k}. 
$$
\ee   
\end{theorem}

\begin{proof}
Since~\eqref{it: hAform decomposition} is a direct consequence of Lemma~\ref{lem: hAform}~\eqref{it: hAform 5}, we focus on~\eqref{it: hAform sym nond}.

We first deal with the symmetric property. Let $x,y$ be  strongly 
homogeneous elements in $\hcalA$. 
We shall prove 
$$ 
\Mn(x\ocalD (y) )=\Mn(y\ocalD (x))
$$
by induction on $\het(y)$. It is obvious when $\het(y)=0$, so we assume that $\het(y)>0$. Since $y$ is  strongly homogeneous, we may assume further that
$y= y' f_{i,p}$ for some $i\in I$,  $p \in \Z$ and a strongly homogeneous $y'$ with $\het(y')<\het(y)$.  Then we have
\begin{align*}
\Mn(x\ocalD(y)) & = \Mn(x\ocalD(y')f_{i,p+1}) = \Mn(f_{i,p-1}x\ocalD(y')) \allowdisplaybreaks\\
& = \Mn(y'\ocalD(f_{i,p-1}x)) = \Mn(y'f_{i,p}\ocalD(x)) = \Mn(y\ocalD(x)), 
\end{align*}
where the second identity follows from Lemma~\ref{Lem: Mfx = Mxf}, the third identity follows from the induction hypothesis. 

For non-degeneracy, it suffices to show that $\hAform{ \ , \ }$ is non-degenerate on each $\hcalA[m]$ by~\eqref{it: hAform decomposition}.  
We shall prove the following claim by induction on $\het(x)$:
\begin{align} \label{eq: non-d claim}
\text{ if $x \in \hcalA[m]$ such that $\hAform{y,x}=0$ for all $y \in \hcalA[m]$, then $x=0$. }    
\end{align}
The claim is obvious when $\het(x)=0$. We assume that
 $\het(x)>0$  and \eqref{eq: non-d claim} holds for $x' \in \hcalA[m]$ with $\het(x')< \het (x)$. 

Assume that $\hAform{y,x}=0$ for all $y \in \hcalA[m]$. 
Then,  Lemma~\ref{lem: hAform}~\eqref{lem: properties of ( , ) (iv)} says that
$$
\hAform{y,\rmE_{i,m}(x)} = \hAform{f_{i,m}y,x}=0 \quad \text{ for any } y \in \hcalA[m]\text{ and } i \in I. 
$$
Lemma~\ref{lem: E Es}~\eqref{it: E and bfL} tells us that $\het(\rmE_{i,m}(x)) < \het(x)$. By the induction hypothesis, we obtain
$\rmE_{i,m}(x)=0$ for any $i\in I$. Then Lemma~\ref{lem: E Es}~\eqref{it: E and bfL} says that $x \in \Qh$. As $\het(x)=0$, we conclude that $x=0$ as we desired. 
\end{proof}

We now define  another symmetric bilinear form on $\hcalA$ as follows:
\begin{align} \label{Eq: def of pair}
\pair{x,y} \seteq \hAform{\sigma(x),\sigma(y)} = q^{-N(\wt(x))}\hAform{x,y} \ \  \text{ for any homogeneous } x,y \in \hcalA,    
\end{align}
where $N$ and $\sigma$ are given in~\eqref{eq: N} and~\eqref{eq: sigma} respectively. Note that 
\begin{align*}
\hAform{f_{i,p},f_{i,p}}=1-q_i^2 \qtq \pair{f_{i,p},f_{i,p}} = q_i^{-1}-q_i,    
\end{align*}
and more generally
\begin{align*}
\hAform{f_{i,p}^n,f_{i,p}^n} =\prod_{k=1}^n (1-q_i^{2k}) \qtq \pair{f_{i,p}^n,f_{i,p}^n} = q_i^{-n^2} \prod_{k=1}^n (1-q_i^{2k}).     
\end{align*}

\begin{lemma} \label{lem: pair}
  Let $i\in I$ and $m\in \Z$ and $n\in\Zp$.
For any $x,y \in \hcalA_{\le m}$, we have 
$$
\pair{f_{i,m}^{(n)}x,y} = q_i^{-n^2} q^{n(\al_{i,m},\wt(x))}\pair{x,\rmE_{i,m}^{(n)}(y)}= q_i^{n^2} q^{n(\al_{i,m},\wt(y))}\pair{x,\rmE_{i,m}^{(n)}(y)}. 
$$

\end{lemma}

\begin{proof}
  By Lemma~\ref{lem: hAform}~\eqref{lem: properties of ( , ) (iv)}, we have
\begin{align*}
\pair{f_{i,m}^n x,y} & = q^{-N(\wt(y)}\hAform{f_{i,m}^nx,y} = q^{-N(\wt(y)}\hAform{x,\rmE_{i,m}^n y} \allowdisplaybreaks\\
& = q^{-N(\wt(y))+N(\wt(x))}\pair{x,\rmE_{i,m}^n y}.
\end{align*}
Since we may assume that $\wt(x)=\wt(y)+n\al_{i,m}$,
we have the desired result.

\end{proof}

\begin{proposition} \hfill \label{prop: pair}
\bnum
\item \label{it: pair (i)}
The bilinear form $\pair{\ , \ }$ is symmetric and non-degenerate.
\item \label{it: pair (ii)}
Let 
 $x=\oprod_{k\in[a,b]}x_k$ and
$y=\oprod_{k\in[a,b]}y_k$ with $x_k,y_k\in\hcalA[k]$,
we have
$$
\pair{x,y} = \prod_{k\in[a,b]} \pair{x_k,y_k}.
$$
\item \label{it: pair (iii)} 
For any $x,y\in \Aqn$ and $m\in\Z$, we have 
$$
\Aform{x,y} = \pair{\varphi_m(x), \varphi_m(y) }
$$
\ee
\end{proposition}

\begin{proof}
The assertions~\eqref{it: pair (i)} and~\eqref{it: pair (ii)} follows from Theorem~\ref{thm: hAform} and~\eqref{eq: N-properties}.     
We shall use induction on $\het(x)$ for proving~\eqref{it: pair (iii)}.
As it is obvious when $\het(x)=0$, we assume that $\het(x) >0$.  

 To prove~\eqref{it: pair (iii)}, it suffices to show that
\begin{align} \label{eq: ets}
\pair{\vph_m(\ang{i}x),\vph_m(y)} = \Aform{\ang{i}x,y}.     
\end{align}
We set
\eqn
&&  \xc \seteq \iota^{-1}(x) \quad   X\seteq\vphi_m(x),\quad\yc \seteq \iota^{-1}(y), \quad Y  \seteq \vph_m(y), \quad
\beta=\wt(x),\quad \gamma \seteq \wt(y).
\eneqn
Note that $\yc \in \Uqgm$ and $X,Y \in \hcalA[m]$. Since the both
sides in~\eqref{eq: ets} are zero unless $\ga=\be -\al_i$, 
we may assume that $\ga =\be-\al_i$. It follows from Lemma~\ref{lem: zeta} that
\begin{align*}
  \aform{\iota(\xc),\iota(e'_i\yc)} &= \zeta^{\be}\kform{\xc,e'_i\yc}
=\zeta^{\be}\kform{f_i\xc,\yc} = \zeta_i\aform{\iota(f_i\xc),\iota(\yc)}  \allowdisplaybreaks\\
& = \aform{\ang{i}\iota(\xc),\iota(\yc)}.
\end{align*}
Thus we have
\begin{align*}
\pair{\vph_m(\ang{i}x),\vph_m(y)} &= q_i^{1/2}\pair{f_{i,m}X,Y} = q_i^{-1/2+\lan h_i,\be \ran} \pair{X,E_{i,m}Y} \allowdisplaybreaks\\
& = \pair{X,\vph_m(\iota(e'_i\yc))} = \aform{ x,\iota(e'_i\yc) } \allowdisplaybreaks\\
& =\aform{\ang{i} x,\iota(\yc)} = \aform{\ang{i}x,y}, 
\end{align*}
 where the  second identity follows from Lemma \ref{lem: pair},
the third from Lemma~\ref{lem: E Es}
and the forth from the induction hypothesis. 
\end{proof}

\begin{remark}
  Proposition \ref{prop: pair} says that the bilinear form $\pair{\ , \ }$ coincides with the bilinear form introduced in \cite[Definition 4.5]{OP24}
 when $ \sfC$ is of finite type. 
\end{remark}

\vskip 2em

\section{Global basis for $\hcalA$} \label{Sec: GB}

Recall that, for each $k \in \Z$, we have  the $\Qh$-algebra isomorphism defined in~\eqref{eq: vph}
$$
\varphi_{k} \cl \Aqhn   \isoto  \hcalA[k] \qquad \text{ by } \varphi_k(\ang{i}) = q_i^{1/2} f_{i,k}.  
$$

One can easily prove the following lemma.

\begin{lemma} \label{lem: vph_k}
The homomorphism $\vph_k$ satisfies the following properties:
\bna
\item $\wt(\vph_k(x)) = (-1)^k \wt(x)$ for any homogeneous element $x \in \Aqn$,
\item \label{it: c auto} we have
\begin{align*}
\varphi_k \circ \sfc(x) = c \circ \vph_k(x) \quad\text{ for any }x\in \Aqn.     
\end{align*}
\ee 
\end{lemma}

We define
$$
\hcalA[k]_\Zq \seteq \vph_{k}(\Azn) \subset \hcalA,
$$
where $\Azn$ is the $\Zq$-lattice defined in  \eqref{eq: Aqn integral form}, and set
$$
\hcalA[a,b]_\Zq \seteq \oprod_{k \in [a,b]} \hcalA[k]_\Zq\subset \hcalA,  \qquad 
\hcalA_\Zq \seteq \sbcup_{a \le b} \hcalA[a,b]_{\Zq}\subset \hcalA. 
$$

\begin{proposition} \label{prop: hcalA integral form}\hfill
\bnum
\item \label{it: Az lattice}
$\hcalA_\Zq$ is a $\Zq$-subalgebra of $\hcalA$, and $\Qh \tens_\Zq \hcalA_\Zq \isoto \hcalA$. 
\item \label{it: c auto on Az}
$\hcalA_\Zq$ is invariant under the action of $c$.
\ee    
\end{proposition}

\begin{proof}
\eqref{it: Az lattice} In order to prove the assertion, it is enough to show that $\hcalA_\Zq$ is closed under the multiplication.

Note that
$$\hcalA[p]_{\Zq}\hcalA[m]_{\Zq}\subset\hcalA[m]_{\Zq}\hcalA[p]_{\Zq}\quad \text{for $p<m-1$}$$
since $xy=q^{-(\wt (x),\wt (y))}yx$ for $x\in\hcalA[p]$ and $y\in\hcalA[m]$.
  
Hence, by Corollary~\ref{cor: main1} and the defining relation~\eqref{Eq: def of hA}~\eqref{Eq: def of hA (b)} of $\hcalA$, it suffices to prove that
$$
\hcalA[m-1]_\Zq \hcalA[m]_\Zq \subset \hcalA[m]_\Zq\hcalA[m-1]_\Zq \quad \text{ for any } m\in \Z. 
$$
Note that~\eqref{eq: aform} and Proposition~\ref{prop: pair}~\eqref{it: pair (iii)} say that
$$
\pair{\vph_k(a),\vph_k(\iota(u_\circ))} = \aform{a,\iota(u_\circ)} = \bAng{a, u_\circ} \quad\text{ for any } a \in \Aqn \text{ and } u_\circ \in \Uqgm. 
$$
Thus we have
\begin{align*}
& \hAmz{m}\hAmz{m-1} \\
& \qquad   = \left\{  x \in \hAm{m-1,m} \ | \  \pair{x,uv} \in \Zq  \right. \allowdisplaybreaks \\
& \hspace{15ex} \left. \text{ for any } u \in \vph_m \circ \iota\big(\Uzgm\big) \text{ and } v \in \vph_{m-1} \circ \iota\big(\Uzgm\big) \right\}.
\end{align*}

Let $x \in \hAmz{m-1}$ and $y \in \hAmz{m}$. We shall prove 
\begin{align} \label{eq: wts}
\pair{xy,uv} \in \Zq    
\end{align}
for any $u \in \vph_m \circ \iota\big(\Uzgm\big)$ and $v \in \vph_{m-1} \circ \iota\big(\Uzgm\big)$ by using induction on $\het(u)$; that is,
we will show that~\eqref{eq: wts} holds for
$$
u = \vph_m \circ \iota(f_i^{(n)}u'_\circ) \qt{with $u_\circ' \in \Uzgm$}
$$ 
under the assumption that~\eqref{eq: wts} holds when $u$
replaced by $\vph_m \circ \iota(u'_\circ)$.

Let $y_\circ \in \Uqgm$ such that $y = \vph_m \circ \iota(\yc)$. Since $y \in \hAmz{m}$, Lemma~\ref{Lem: mem for An and ei'} and Lemma~\ref{lem: E Es}~\eqref{it: En}
say that
\begin{align} \label{eq: E(n) and Anz}
q_i^{n/2}\zeta_i^{-n} \rmE_{i,m}^{(n)}(y) \equiv_q \vph_m ( \iota( \zeta_i^{-n} e_i^{\prime(n)}\yc  )  ) \in \hAmz{m}.    
\end{align}

Setting $u' \seteq \vph_m \circ \iota(u'_\circ) \in \vph_m \circ \iota\big(\Uzgm\big)$, we have
$$
u = \vph_m \circ \iota(f_i^{(n)}u_\circ') = q_i^{n/2} \zeta_i^{-n} f_{i,m}^{(n)} u'.
$$
Since $\rmE_{i,m}(x)=0$ by Lemma~\ref{lem: E Es}~\eqref{it: E comm}, 
we have
\begin{align*}  
\pair{xy,uv} &\equiv_q q_i^{n/2} \zeta_i^{-n} \pair{xy,f_{i,m}^{(n)} u'v} \allowdisplaybreaks\\
&\equiv_q q_i^{n/2} \zeta_i^{-n} \pair{\rmE_{i,m}^{(n)}( xy), u'v} \allowdisplaybreaks\\
&\equiv_q \pair{x\rmE_{i,m}^{(n)}( q_i^{n/2} \zeta_i^{-n}  y), u'v} \in \Zq
\end{align*}
where the second identity follows from Lemma \ref{lem: pair},
the third identity follows from Lemma \ref{lem: E Es}~\eqref{it: E Es} and the last one follows from \eqref{eq: E(n) and Anz}.

\smallskip
\noindent
\eqref{it: c auto on Az} follows from Lemma~\ref{lem: vph_k}~\eqref{it: c auto}.
\end{proof}

\begin{remark}
Note that the product $\oprod \vph_k \circ \iota(\Uzgm) \subset \hcalA$ is not a $\Zq$-subalgebra of $\hcalA$, as seen by 
$$
\vph_k(\iota(f_i))\vph_{k+1}(\iota(f_i)) = q_i^2 \vph_{k+1}(\iota(f_i)) \vph_k(\iota(f_i)) +\dfrac{1}{q_i^{-1}-q_i}.
$$
\end{remark}

We now define $\Z[q]$-lattices as follows:
\begin{equation} \label{Eq: def of Zq lattices of hA}
\begin{aligned}
\Lup\left( \hAmz{k} \right) & =\vph_k\left(\Lup\big(\Azn\big)\right),  \allowdisplaybreaks\\
\Lup\left( \hAmz{a,b} \right) & = \oprod_{k \in [a,b]} \Lup\left( \hAmz{k} \right),  \allowdisplaybreaks\\
\Lup\left( \hcalA_\Zq \right) & = \sbcup_{a \le b} \Lup\left( \hAmz{a,b} \right). 
\end{aligned}
\end{equation}
Recall the \emph{extended crystal} of the infinite crystal $B(\infty)$
introduced in \S\;\ref{subsec: crystals}:
$$
\hB(\infty) =\{ (b_k)_{k\in \Z} \ | \  b_k \in B(\infty) \text{ and } b_k = \mathsf{1} \text{ except for finitely many } k  \},
$$
where $\mathsf{1}$ is the highest weight vector in $B(\infty)$ of $\Uqgm$. For any $\bfb =(b_k)_{k\in\Z} \in \hB(\infty)$, we set
$$
\rmP(\bfb) \seteq \oprod_{k \in \Z} \vph_k(\Gup(b_k)) \in \Lup(\hcalA_\Zq). 
$$
Then, $\st{\rmP(\bfb)\mid\bfb\in \hBi}$ is a $\Z[q]$-basis of $\LuphA$.

We regard $\hB(\infty)$ as a $\Z$-basis of $\Lup(\hcalA_\Zq)/q\Lup(\hcalA_\Zq)$ by 
$$
\bfb  = \rmP(\bfb) \ {\rm mod} \ q\Lup(\hcalA_\Zq). 
$$

\begin{lemma} \label{lem: [x,y]q}
For any $x \in \hAm{k}_\al$ and $y \in \hAm{k+1}_\be$, we have
$$
[x,y]_q  \in \sum_{\substack{ \al',\be' \\ \dVert{\al'} < \dVert{\al}, \; \dVert{\be'} < \dVert{\be} }} \hAm{k+1}_{\be'}\hAm{k}_{\al'},
$$
where $\dVert{ \cdot }$ is given in~\eqref{eq: het dVert}, and the order $<$ on $\rl$ is given in~\eqref{eq: < order on Q}. 
\end{lemma}

\begin{proof}
We shall use induction on $\het(y)$. 
Since it is obvious when $\het(y)=0$, we may assume that $\het(y)>0$. We may assume that $y = f_{i,k+1} y'$ for $y'\in \hAm{k+1}$. 
It follows from the defining relation \eqref{Eq: def of hA} \eqref{Eq: def of hA (b)} of $\hcalA$ that  
\begin{align*}
[x ,f_{i,k+1}]_q  = \sum_{\dVert{\al'}  <  \dVert{\al} } \hAm{k}_{\al'}, 
\end{align*}	
which says that, by the induction hypothesis and~\eqref{eq: []q1},  
\[
[x,y]_q  = [x, f_{i,k+1} ]_q \;  y'     + 	q^{ - ( \al, -\al_{i,k+1} )}  f_{i,k+1} [x,y']_q  
	 \in \hspace{-5ex} \sum_{\substack{ \al',\beta'\in\rl\\ \dVert{\al'}<\dVert{\al}, \;  \dVert{\be'}<\dVert{\be}}} \hspace{-5ex} 
	\hcalA[k+1]_{\beta'}\hcalA[k]_{\al'}. \qedhere 
\]	    
\end{proof}

For $\bfb \in (b_k)_{k\in \Z} \in \hB(\infty)$ and  $\bfb' \in (b'_k)_{k\in \Z} \in \hB(\infty)$, we define 
\begin{align} \label{Eq: order on hB}
  \bfb\preceq\bfb' \qt{$\dVert{\wt(b_k)}\le\dVert{\wt(b'_k)}$ for any $k \in \Z$.}
\end{align}	

Note that $\preceq$ in \eqref{Eq: order on hB} is a preorder on $\hB(\infty)$,
and we introduce the relation $\prec$ as in {\bf Convention}~\eqref{conv:preorder}.

\begin{lemma} \label{lem: PP}
For any $\bfb =(b_k)_{k\in \Z} \in \hBi$, we have
$$
c(\rmP(\bfb)) - \rmP(\bfb) = \sum_{\bfb' \in \hBi, \; \bfb' \prec \bfb} a_{\bfb,\bfb'}(q) P(\bfb')
$$
for some $a_{\bfb,\bfb'}(q) \in \Zq$. 
\end{lemma}

\begin{proof}
Let $\hBel \in \hBi$. Then there exist integers $a\le b$ such that $b_k =\mathsf{1}$ unless $k \in [a,b]$. Thus we write
$$
\rmP(\bfb) = G_bG_{b-1} \cdots G_a,
$$
where we set $G_k \seteq \vph_k(\Gup(b_k))$ for $k \in [a,b]$. Let $\be_k \seteq \wt(G_k)\in\rl$. It follows from~\eqref{eq: c anti} and Lemma~\ref{lem: vph_k} that
\begin{align*}
c(\rmP(\bfb)) & = q^{\sum_{a<b} (\be_a,\be_b)} c(G_a)\cdots c(G_{b-1})c(G_b) \allowdisplaybreaks \\
&  = q^{\sum_{a<b} (\be_a,\be_b)}  G_a \cdots  G_{b-1} G_b.
\end{align*}
Therefore, the assertion follows from Proposition \ref{prop: hcalA integral form} and Lemma \ref{lem: [x,y]q}.
\end{proof}

The following theorem is a direct consequence of Lemma \ref{lem: PP}. For the reader's convenience, we write a short proof of the theorem (see \cite[Section 5]{Tingley17} and references therein). 

\begin{theorem} \label{Thm: global basis}\ 
\bnum
\item \label{it: global (i)}
For each $\hBel \in\hBi$, there exists a unique
$\rmG(\bfb)\in \LuphA$ such that
\bna
\item
  $\rmG(\bfb)-\rmP(\bfb)\in \displaystyle\sum_{\bfb'\prec\bfb}q\Z[q]\rmP(\bfb')$,
\item
  $c(\rmG(\bfb))=\rmG(\bfb)$.  
  \ee
\item \label{it: global (ii)}  Moreover, the set $\{\rmG(\bfb) \ | \ \bfb\in\hBi\}$ forms a $\Z[q]$-basis of $\LuphA$, and
 a $\Z$-basis of $\LuphA\cap c\big(\LuphA\big)$. 
 \item \label{it: global (iii)}
 For any $\bfb \in\hBi$, we have 
 $$
 \rmP(\bfb) = \rmG(\bfb) + \sum_{\bfb'  \prec \bfb}  f_{\bfb,\bfb'}(q) \rmG(\bfb')\qquad \text{ for some $f_{\bfb,\bfb'}(q) \in q \Z[q]$.}
 $$
  \ee
\end{theorem}

\begin{proof}
\eqref{it: global (i)} 
Note that, for each $\hBel \in \hBi$ there exist  $r\le s \in \Z$ 
such that $b_k =\mathsf{1}$ unless  $k \in [r,s]$  and $b_r,b_s \ne \mathsf{1}$.
We shall prove by using induction on $s-r$.     

We first check the case $s-r=0$; i.e., $b_k =\mathsf{1}$ for all $k \in \Z \setminus \{s\}$.  Then we have
$\rmG(\bfb) = \rmP(\bfb)=\vph_s(\rmG(b_s))$. Then it satisfies {\rm (a)} and {\rm (b)} trivially. 

We now consider a general case. By Lemma~\ref{lem: PP} and the induction hypothesis, we write
$$
c(\rmP(\bfb)) = \rmP(\bfb) + \sum_{\bfb' \prec \bfb} p_{\bfb,\bfb'}(q)\rmG(\bfb')
$$
for some $p_{\bfb,\bfb'}(q) \in \Zq$. Since $c$ is an involution,
$p_{\bfb,\bfb'}(q^{-1})=-p_{\bfb,\bfb'}(q)$, and hence
$$
p_{\bfb, \bfb'}(q) = f_{\bfb, \bfb'}(q) - f_{\bfb, \bfb'}(q^{-1})
$$
for some polynomial $f_{\bfb, \bfb'}(q) \in q\Z[q]$.
We define
$$
\rmG(\bfb) \seteq \rmP(\bfb) + \sum_{\bfb'\prec \bfb} f_{\bfb,\bfb'}(q) \rmG(\bfb'). 
$$
One can easily check that $c(\rmG(\bfb)) = \rmG(\bfb)$, which tells us that $\rmG(\bfb)$ satisfies the desired properties. 

\smallskip
\noindent
\eqref{it: global (ii)} Corollary \ref{cor: main1} says that  $\{ \rmP(\bfb) \mid \bfb \in \hBi \}$ forms a basis of $\LuphA$. Thus the assertion follows from \eqref{it: global (i)}. 

\smallskip
\noindent
\eqref{it: global (iii)}  It follows directly from~\eqref{it: global (i)} and~\eqref{it: global (ii)}.
\end{proof}

Thanks to Theorem \ref{Thm: global basis}, the \emph{global basis} of $\hcalA$ is defined as follows.

\begin{definition} \ \label{def: global basis}
  \bnum
\item   We call  $\{\rmG(\bfb)   \mid   \bfb\in\hBi \}$ the \emph{global basis} of $\hcalA$.
  \item  We define
  $$
  \trmG(\bfb) \seteq q^{ -N(\wt( \rmG(\bfb) )/2} \rmG(\bfb) \qquad \text{ for any $\bfb \in \hBi$,}
  $$
  and call $\{\trmG(\bfb)   \mid  \bfb\in\hBi \}$ the \emph{normalized global basis} of $\hcalA$.
  \ee    
\end{definition}
\noindent
Note that the normalized global basis elements $\trmG(\bfb)$ are bar-invariant by \eqref{Eq: bar inv c inv}.

\begin{proposition} 
For any $\bfb \in \hB(\infty)$, we have 
$$
\rmG(\bfb)^\star = \rmG(\bfb^\star) \qquad \text{ and }\qquad \trmG(\bfb)^\star = \trmG(\bfb^\star).
$$	
\end{proposition}
\begin{proof}
By the definition \eqref{eq: star} of $\star$, we have 
$$
 \varphi_{k} (x)^\star = \varphi_{-k} (x ^*)  \qquad \text{ for any $x \in \Aqhn$ and $k\in \Z$.} 
$$	
Hence, comparing the action of $\star$ on $\hBi$ in \eqref{def:starc},  we obtain
$$
\rmP(\bfb) ^ \star = \rmP(\bfb^ \star), 
$$
which completes the proof by Theorem \ref{Thm: global basis}.
\end{proof}

Set
$$
\hcalA_{\Zqh}=\Zqh\tens_{\Zq}\hcalA_{\Zq},
$$
and, for any $m\in\Z$ and $b \in B(\infty)$, define $\bfe_m(b)=(b_k)_{k\in\Z} \in \hBi$ by 
\begin{align} \label{Eq: emb} 
b_k \seteq \begin{cases}  b & \text{ if } k=m, \\  \mathsf{1} & \text{ otherwise.}
\end{cases}	
\end{align}

The following proposition is now obvious. 

\begin{proposition} \label{Prop: properties of Gb} \
  \bnum
  \item \label{Prop: properties of Gb (i)}
  For any $m\in \Z$ and $b \in B(\infty)$, we have $ \vph_m( \Gup(b)) = \rmG(\bfe_m(b))$. 
  \item
The set $\{\trmG( \bfb)\mid \bfb\in\hBi\}$ is a $\Zqh$-basis of $\hcalA_{\Zqh}$.
\item
For $ \bfb, \bfb'\in\hBi$, we have
$$\hAform{ \trmG(\bfb),\trmG(\bfb') } \in\Z[[q]]\cap\Q(q).$$
\item For $\bfb, \bfb'\in\hB(\infty)$, we have $\hAform{ \trmG(\bfb),\trmG(\bfb')} \big\vert_{q=0}=\delta_{\bfb, \bfb'}$,
\item We have
\begin{align*}
& \{ x \in \hcalA_\Zqh \mid \ol{x} = x, \ \hAform{x,x} \in 1 + q^{1/2}\Q[[q^{1/2}]] \}   \allowdisplaybreaks\\
& \hspace{8ex} = \{ \trmG(\bfb) \ | \ \bfb \in \hBi\} \cup \{ -\trmG(\bfb) \ | \ \bfb \in \hBi\}. 
\end{align*}
\ee 
\end{proposition}

\vskip 2em 

\section{Quantum Grothendieck ring and global bases} \label{Sec: QGR and GB}

\subsection{Quantum Grothendieck ring} \label{Sec: qGr}
In this subsection, we briefly recall the quantum Grothendieck ring of the
skeleton category of the quantum affine algebra, which is also referred to as the \emph{Hernandez-Leclerc} category. 
We basically follows the notations given in \cite[Section 6.1]{KKOP24}. 

Let $t$ be an indeterminate and let $U_t'(\g)$ be the quantum affine algebra whose bases field $\bfk$ is the algebraic closure
of $\Q(t)$  in $\cup_{m>0} \C\llf t^{1/m}\rrf$.\footnote{By the technical reason, 
we swap the roles of $t$ and $q$ in the common notations related to the quantum Grothendieck ring over quantum affine algebras (see \cite{FHOO22} for instance).}
We denote by  
$\Cg$ the category of finite-dimensional integrable modules over  $U_t'(\g)$.
We set $I_0 \seteq I \setminus \{ 0 \}$, where we refer the reader to  \cite[Section 2.1]{KKOP24} (see also \cite[Section 2.3]{KKOP23P}) for the choice of $0$. We denote by $\g_0$ the Lie subalgebra of $\g$ generated by $f_i$, $e_i$ and $h_i$ for $i\in I_0$. For any module $M \in \Cg$, we write $\rdual(M)$ for the \emph{right dual} of $M$ and set $\rdual^k(M)\seteq  \rdual( \rdual^{k-1}(M)) $ for $k\in \Z$. For each $i\in I_0$, let $\varpi_i$ be the 
\emph{$i$-th fundamental weight} and denote by $V(\varpi_i)_a$
the \emph{$i$-th fundamental weight module}  with \emph{spectral parameter} 
$a \in \bfk^\times$.

Let $\gf$ be the simple Lie algebra of simply-laced type associated with the quantum affine algebra $U_t'(\g)$ (see \cite{KKOP22} and \cite[Section 6.1]{KKOP24}), and let $\sfC_{{\fin}}$ be the corresponding Cartan matrix.
 Let $I_\fin$ be the index set for simple roots of $\g_\fin$. 
We choose  a $Q$-datum  $\calQ = \Qdatum$ associated to $\g_0$, which consists of the \emph{Dynkin diagram} $\Dynkin$ of $\gf$, the \emph{folding automorphism} $\varsigma$ between $\gf$ and $\g_0$  and a \emph{height function} $\xi\cl I_{{\rm fin}} \rightarrow \Z$ (see \cite{FO21} and \cite[Section 6.1]{KKOP24} for precise definitions). 
For vertices $\im, \jm \in \Dynkin$, we denote by $d(\im, \jm)$ the distance between $\im$ and $\jm$ in the Dynkin diagram $\Dynkin$.
The automorphism $\varsigma$ gives a surjective map 
$$
\pi\cl I_{{\rm fin}} \twoheadrightarrow I_0.
$$
Note that every fiber of $\pi$ is a $\varsigma$-orbit in $I_\fin$. 
We set 
$$
\sigma_0(\g) \seteq \{  (\imath, p) \in I_{{\rm fin} } \times \Z \mid p-\xi(\imath) \in 2 d_\imath \Z  \},
$$
where $d_\imath = | \pi^{-1} (\pi(\imath))| $. 

Recall that the isomorphism classes in $\Cg$ are parameterized by the set $(1+z\bfk[z])^{I_0}$ of $I_0$-tuples of monic polynomials, called Drinfeld polynomials \cite{CP95A, CP95},
where $z$ is an indeterminate.  Let us introduce a formal variable $Y_{\pi(\im),p}$ for each $(\im,p) \in \sigma_0(\g)$ and consider the Laurent polynomial ring $\calY \seteq \Z[Y_{\pi(\im),p}^{\pm 1} \mid (\im,p) \in \sigma_0(\g)]$. Let $\calM \subset\calY$ be the set of all Laurent monomials. 
We write $m \in \calM$ as 
\begin{align*}
m = \prod_{(\im,p) \in \sigma_0(\g)} Y_{\pi(\im),p}^{u_{\pi(\im),p}(m)}.
\end{align*}
We say that an element $m \in \calM$ is \emph{dominant} if $u_{\pi(\im),p}(m) \ge 0$ for all $(\im,p) \in \sigma_0(\g)$ and set $\calM^+ \subset \calM$ the set of all dominant monomials.
For each $m \in \calM^+$, we have a simple module $L(m) \in \Cg$ corresponding to the Drinfeld polynomial $(\prod_p (1-t^pz)^{u_{\pi(\im),p}(m)})_{\pi(\im) \in I_0}$.
Note that $L(Y_{i,p})$ corresponding to the dominant monomial $Y_{i,p}$ is a fundamental representation.

Note that we have the following bijection determined by the \emph{$\calQ$-adapted Coxeter element}:
$$
\phi_\calQ\cl \sigma_0(\g) \buildrel \sim \over \longrightarrow \Phi^+_\gf \times \Z, 
$$
where $\Phi^+_\gf$ is the set of positive roots of $\gf$ (see \cite{FO21}).
Setting $ \sigma_Q(\g)  \seteq \phi_\calQ^{-1} \left( \Phi^+_\gf\times \{0\} \right)$, we have the bijection 
\begin{align} \label{Eq: 1-1 Detla Q}
\phi_\calQ|_{\sigma_Q(\g)} \cl 	\sigma_Q(\g) \buildrel \sim \over \longrightarrow \Phi^+_\gf,
\end{align}
where we identify $\Phi^+_\gf$ with $\Phi^+_\gf\times \{0\}$.
For each $(\imath, p) \in \sigma_0(\g)$ and $\beta \in \Phi^+_\gf$, we define 
$$V(\imath, p) \seteq L(Y_{\pi(\im),p}) \qtq V_\calQ(\beta) \seteq V(p_\beta)$$
where $p_{\beta} \in \sigma_0(\g)$ denotes the pair corresponding to $\beta$ in the bijection \eqref{Eq: 1-1 Detla Q}. 

The \emph{Hernandez-Leclerc category} $\Cg^0$ (resp.\ $\scrC_\calQ$) is the smallest full subcategory of $\Cg$ that 
\bna
\item contains $V(\imath, p)$ for all $(\imath, p) \in \sigma_0(\g)$ (resp.\ $(\imath, p) \in \sigma_\calQ(\g)$), and
\item is stable by taking subquotients, extensions and tensor products. 
\ee
Note that the set $\{  \rdual^{m} V_\calQ(\beta) \mid m\in \Z, \beta\in\Phi^+_\gf \}$ 
coincides with the set of all fundamental modules in $\Cg^0$.

\smallskip 
We now assume that $U_t'(\g)$ is of untwisted affine type.
In \cite{FR99}, Frenkel-Reshetikhin constructed an injective ring homomorphism 
$$
\chi_t \cl K(\Cg^0) \to \calY
$$
which is known as the \emph{$t$-character homomorphism}. Here $K(\Cg^0)$ denotes the Grothendieck ring of $\Cg^0$. As a ring $K(\Cg^0)$ is isomorphic to the
ring of the polynomials in $\{ [L(Y_{\pi(\im),p})] \}$ (\cite{FR99}). 

Let $\calY_q$ be the the \emph{quantum torus} associated with $\g_0$,
which is the $\Zqh$-algebra
generated by $\st{\tY_{\pi(\im),a}^{\pm 1} \mid(\im,a) \in \sigma_0(\g)}$
with the following defining relations:
\begin{align*}
\tY_{\pi(\im),p}\tY_{\pi(\im),p}^{-1} & = \tY_{\pi(\im),p}^{-1}\tY_{\pi(\im),p}=1, \\
\tY_{\pi(\im),p} \tY_{\pi(\jm),s} &= q^{ k+l+\delta(p<s) } \delta((\pi(\im),p)\ne (\pi(\jm),s)) \tY_{\pi(\jm),s}\tY_{\pi(\im),p},    
\end{align*}
where $\phi_\calQ^{-1}(\im,p)=(\al,k)$ and  $\phi_\calQ^{-1}(\im,p)=(\be,l)$\footnote{The definition of quantum torus in \cite{HL15,FHOO22} is different from ours and can be obtained by replacing
$q^{\pm 1/2}$ with $t^{\mp 1/2}$.}
(see \cite{Her03, Her04, Nak04, VV03} and see also \cite[Section 5]{HL15} for details).

Note that 
the quantum torus $\calY_q$ has a \emph{bar-involution} fixing $\tY_{\pi(\im),p}$ and sending $q^{\pm1/2}$ to $q^{\mp1/2}$,

For a simple module $L \in \Cg^0$, we denote by   $\chi_{t,q}(L)$  the \emph{$(t,q)$-character} or the \emph{$q$-deformed $t$-character} of $L$.
The $(t,q)$-character $\chi_{t,q}(L)$ of $L$ is contained in $\calY_q$ 
and \emph{bar-invariant}.

We write  $\calK_q(\Cg^0)$  (resp.\ $\calK_q( \rdual^m (\scrC_\calQ))$) for the \emph{$q$-deformed} (or \emph{quantum}) \emph{Grothendieck ring} of $ \Cg^0$ 
(resp.\ $\rdual^m (\scrC_\calQ)$) over $\Zqh$. 
They are the $\Zqh$-subalgebra of $\calY_q$ generated by
$\chi_{t,q}(L)$ for simple $L$'s in $\Cg^0$ (resp.\ $ \rdual^m (\scrC_\calQ))$).

We now define
\begin{equation} \label{Eq: K and varkappa}
\begin{aligned}
\bbK_{m}^\calQ &\seteq     \Qh \otimes_{\Zqh }  \calK_q( \rdual^m (\scrC_\calQ)),   \\
\bbK^0 &\seteq   \Qh \otimes_{\Zqh } \calK_q(\Cg^0).
\end{aligned}		
\end{equation} 
When $m=0$, we simply write $ \bbK^\calQ = \bbK_{0}^\calQ$.

\subsection{Global basis and simple modules} \label{Sec: gb sm}
Let $\Aqn$ be the quantum unipotent coordinate ring associated with the Cartan matrix $\sfC_{\fin}$ and let $\hcalAf$ be the bosonic extension of $\Aqn$. 
We fix a $Q$-datum $\calQ$ throughout this subsection. 
The following lemma was proved in  \cite[Theorem 6.1]{HL15}, \cite[Corollary 10.2]{HO19} and  \cite[Corollary 8.7]{FHOO22}.

\begin{lemma} [\cite{FHOO22, HL15, HO19}] \label{Lem: Upsilon} \
\bnum
\item There is an isomorphism 
 \begin{align*}
 \Upsilon_{\calQ, m} \cl \Aqhn \buildrel \sim \over \longrightarrow \bbK_{ m}^\calQ, \qquad \ang{ \imath } \mapsto   q^{1/2}  \chi_{t,q}( \rdual^m V_\calQ(\al_{\imath})).
 \end{align*}	
\item \label{Lem: Upsilon (ii)}
  The isomorphism $\Upsilon_{\calQ, m}$ takes the \emph{normalized} upper global basis to the set of $t$-deformed $q$-characters of simple modules in $\scrC_\calQ$. More precisely, we have 
$$ 
\Upsilon_{\calQ, m} \left( q^{ -N(\wt(b))/2} \Gup(b) \right) = \chi_{t,q}( \rdual^m V_\calQ(b)) \quad \text{ for any  } b\in B(\infty),
$$ 
where $V_\calQ(b)$ is the simple module in $\scrC_\calQ$ corresponding to $b$ via $\Upsilon_{\calQ, 0 }$. 
\ee
\end{lemma}

\begin{remark} \
\bnum	
\item	The indeterminate $v$ is used in  the quantum unipotent coordinate ring in \cite[Theorem 6.1]{HL15}. The indeterminate $v$ corresponds to $q^{-1}$ in our paper setting.
\item $N(\beta)$ in our paper is different from  $N(\beta)$ used in \cite[(3.14) and Theorem 6.1]{HL15}.
\ee   		
\end{remark}	

The homomorphism  $\Upsilon_{\calQ,m}$ extends to
 an $\Qh$-algebra isomorphism
from the quantum Grothendieck ring to the bosonic extension as follows:
 \begin{align} \label{Eq: OmegaQ}
\Omega_\calQ \cl  \bbK^0  \buildrel \sim \over \longrightarrow \hcalAf, \qquad \chi_{t,q}( \rdual^m V_\calQ(\al_{\imath }) ) \mapsto f_{ \imath ,m}
\end{align}	
(see \cite[Theorem 7.3]{HL15} and \cite[Theorem 10.4]{FHOO22}). 
Note that $\Omega_\calQ$ preserves the bar-involution.
We set 
$$
 \Omega_{\calQ, m} \seteq \Omega_{\calQ}\mid_{ \bbK_{ m}^\calQ} \;\cl\; \bbK_{ m}^\calQ \buildrel \sim \over \longrightarrow \hcalAf[m] \subset \hcalAf.
$$

Then we have the following commutative diagram:  
\begin{equation} \label{Eq: dig for qt and hA}
\begin{aligned} 
\xymatrix@C=10ex{\Aqn \ar[r]^{\Upsilon_{\calQ,m} }  \ar[rd]_{\vphi_m} & \bbK_{ m}^\calQ   \ar@{>->}[r] \ar[d]^ {\Omega_{\calQ, m}} \  & \bbK^0 \ar[d] ^ {\Omega_{\calQ}} \\
	& \hcalAf[m] \  \ar@{>->}[r]  &  \hcalAf
}
\end{aligned}	
\end{equation}
where $\vph_m$ is given in
\eqref{eq: vph}.

The following lemma follows immediately
from Proposition~\ref{Prop: properties of Gb}~\eqref{Prop: properties of Gb (i)}.
\begin{lemma} \label{Lem: qt to tGb}
For any $b \in B(\infty)$, the isomorphism $\Omega_{\calQ, m}$ takes the $(t,q)$-character of $ \rdual^m V_\calQ(b)$ to the normalized global basis element $\trmG( \bfe_m(b))$, i.e., 
$$
\Omega_{\calQ, m} (  \chi_{t,q}( \rdual^m V_\calQ(b)) ) = \trmG( \bfe_m(b)),
$$
where $\bfe_m(b)$ is defined in \eqref{Eq: emb}.
\end{lemma}

Let $\weyl_\fin$ be the Weyl group for $\g_\fin$.
We choose a $\calQ$-adapted reduced sequence $\bfi = (\im_1, \im_2, \ldots, \im_\ell)$ of the longest element $w_0$ in $\weyl_\fin$, and 
denote by $ \fraks \seteq  \left((\im_k, p_k) \right)_{k\in \Z} $ the corresponding \emph{admissible sequence} in $\sigma_0(\g)$ (see \cite[Definition 6.9]{KKOP24} for the definition of an admissible sequence, and see also \cite[Proposition 6.11]{KKOP24}).   It follows from \cite[Proposition 6.11]{KKOP24} that 
\bnum
\item the set $   \{ V(\imath_k, p_k) \mid k\in \Z \} $ coincides with the set of all fundamental modules in $\Cg^0$,
\item the set $  \{ V(\imath_k, p_k) \mid k\in [1,\ell  ] \} $ coincides with the set of all fundamental modules in $ \scrC_\calQ$,
\item $ \imath_{k+\ell} = (\imath_k)^*$ and $ p_{k+\ell} = p_k + \mathrm{ord}(\varsigma) h^\vee$, where $\al_{\im^*} = -w_0(\al_\im)$, $\mathrm{ord}(\varsigma) $ is the order of $\varsigma$ and $h^\vee$ is the \emph{dual Coxeter number} of $\g_0$.
\ee

For any $k\in [1, \ell]$, we set $\beta_k$ be the positive root corresponding to $(\imath_k, p_k)$ via the bijection \eqref{Eq: 1-1 Detla Q}. 
Note that $ \Phi^+_\fin = \{ \beta_\ell >_\bfi \cdots >_\bfi \beta_{2} >_\bfi \beta_1 \}$ and $>_\bfi$ is the convex order on $\Phi^+_\fin$ corresponding to the reduced sequence $\bfi$.

We hereafter regard $B(\infty)$ as a \emph{PBW monomial realization} associated with $\bfi$ so that an element $\bfa \in \Zp^{\oplus [1,\ell]}$ can be understood as an element in $B(\infty)$. Similarly, for a given $\bfc = (c_k)_{k\in \Z} \in \Zp^{\oplus \Z}$, we now  understand $ \bfc_k \seteq ( c_{k \ell+1}, c_{k \ell +2}, \ldots, c_{(k+1)\ell} ) $ as an element of $B(\infty)$ via the PBW monomial realization,  and $\bfc$ as an element of the extended crystal $\hBi$ via the PBW monomial realization (see Section \ref{subsec: crystals}). Thus, we have the following identifications:  
\begin{align} \label{Eq: PBW identiftcation}
B(\infty) \simeq  \Zp^{\oplus[1,\ell]} \qtq \hBi \simeq  \Zp^{\oplus\Z}.
\end{align}  
 
For any $m\in \Z$ and $ k \in [1,\ell]$, we set 
$$
x_m (\beta_k) \seteq  \chi_{t,q}( \rdual^m V_\calQ( \beta_k) ) \in \bbK^0,
$$
and define 
$$ 
F_\bfi(\beta_k) \seteq  q^{ N(\beta_k)/2 } \Upsilon_{\calQ,0}^{-1} ( x_0 (\beta_k) ) \in \Aqn.
$$
Then $ \{   F_\bfi(\beta_k) \mid   k\in [1,\ell] \}$ is the set of dual PBW vectors in $\Aqn$ corresponding to $ \bfi$
(see \cite{KKK2,KO19,OS19}).
For any $\bfa = (a_k) \in \Zp^{\oplus [1,\ell]}$, we define 
$$
M_\bfi(\bfa) \seteq q^A F_\bfi(\beta_\ell)^{a_\ell} \cdots  F_\bfi(\beta_2)^{a_2} F_\bfi(\beta_1)^{a_1},
$$
where $A \in  \Z$ is chosen so that  
\begin{align} \label{Eq: M G}
	M_\bfi(\bfa) \equiv  \Gup_\bfi(\bfa) \qquad \text{mod } q\LupA, 
\end{align}
where we understand $\bfa$ as an element of $B(\infty)$. The monomial $M_\bfi(\bfa)$ is called a \emph{dual PBW monomial}.
The set $\st{M_\bfi(\bfa)\mid\bfa\in \Zp^{\oplus[1,\ell]}}$ is a $\Z[q]$-basis of
$\LupA  $.
We now define 
\begin{align} \label {Eq: pE_m}
\bbE_m(\bfa) \seteq q^{-N(\gamma)/2}   \Upsilon_{\calQ,m} ( M_\bfi(\bfa)),
\end{align}
where $\gamma = \wt( M_\bfi(\bfa))$.
Note that 
$ \bbE_m(\bfa) $ is equal to $x_m (\beta_\ell)^{a_\ell} \cdots x_m (\beta_2)^{a_2} x_m (\beta_1)^{a_1} $ up to a power of $q^{\pm 1/2}$.

For a dominant monomial $x = \prod_{(\im,p) \in \sigma_0(\g)} Y_{\pi(\im),p}^{u_{\pi(\im),p}(x)}$, 
we define 
$$
 \bbE(x) \seteq q^{u_x} \oprod_{p\in\Z} \prod_{(\im,p) \in \sigma_0(\g)} \chi_{t,q}(V(\im, p))^{u_{\pi(\im),p}(x)},
 $$ 
 where $u_x$ is chosen so that $x$ appears with coefficient $1$ in the expansion of $\bbE(x)$ on the basis of \emph{commutative} monomials in the quantum torus $\calY_q$ (see \cite[Section 3.2]{HL15} for the definition of commutative monomials). 
 The element  $\bbE(x)$ is called  the $(t,q)$-character of the \emph{standard module} corresponding to the dominant monomial $x$.

For any $\bfc = (c_k)_{k\in \Z} \in \Zp^{\oplus \Z}$, we denote by $\sfm(\bfc) \seteq \prod_{k\in \Z} Y_{\pi(  \imath_k ), p_k}^{c_k}$ the dominant monomial corresponding to the admissible sequence $\fraks =  \left((\imath_k, p_k) \right)_{k\in \Z}$.   We set $ \bbE(\bfc) \seteq  \bbE(\sfm(\bfc))$. 
By \cite[Proposition 10.1]{FHOO22}, we have 
\begin{align} \label{Eq: standard module}
\bbE(\bfc) =   q^{- \frac{1}{2} \sum_{a>b} (\ga_a, \ga_b) }  \oprod_{k\in\Z} \bbE_k(\bfc_k) \in \bbK^0,
\end{align}
where $ \bfc_k = ( c_{k \ell+1}, c_{k \ell +2}, \ldots, c_{(k+1)\ell} )$ and $\ga_k  = (-1)^k \wt(  M_\bfi(\bfc_k) )  $.

\begin{theorem} \label{Thm: simple to gb}
For any $\bfc = (c_k)_{k\in \Z} \in \Zp^{\oplus \Z}$, we have 
$$
 \Omega_\calQ( \chi_{t,q}( L(\bfc) )) = \trmG(\bfc),
$$
where $ L(\bfc)$ is the simple module in $\Cg^0$ corresponding to  the dominant monomial $\sfm(\bfc)$.
\end{theorem}

\begin{proof}
  For $m\in \Z$ and $\bfa \in \Zp^{\oplus [1,\ell]}$, we set $\bbM_m(\bfa) \seteq  \vph_m( M_\bfi(\bfa))$. 	
Let $\bfc = (c_k)_{k\in \Z} \in \Zp^{\oplus \Z}$ and set  $ \bfc_k \seteq ( c_{k \ell+1}, c_{k \ell +2}, \ldots, c_{(k+1)\ell} )$ for $k\in \Z$. 
Define 
$$
\bbM(\bfc) \seteq \oprod_{k\in\Z} \bbM_k(\bfc_k)\in \hA.
$$
Then
$\st{\bbM(\bfc)\mid \bfc\in\Zp^{\oplus \Z}}$ is a $\Z[q]$-basis of $\LuphA$,
and
we have
$$\bbM(\bfc) \equiv \rmG (\bfc)\ \mod q\LuphA.$$

By \eqref{Eq: pE_m} and \eqref{Eq: standard module}, we have
\eqn
&&\Omega_{ \calQ}(\bbE(\bfc))=q^{-N(\wt( \bfc))/2}\bbM(\bfc),
\eneqn
where $\wt(\bfc)\seteq\wt( \bbM(\bfc))$.
Hence
$\st{q^{N(\wt(\sfc))/2}\Omega_{ \calQ}(\bbE(\bfc))\mid\bfc\in\Zp^{\oplus\Z}}$ is a $\Z[q]$-basis of $\LuphA$
and
$$q^{N(\wt(\sfc))/2}\Omega_{ \calQ}(\bbE(\bfc))\equiv  \rmG  (\bfc)\ \mod q\LuphA.$$

\medskip
On the other hand,  the unitriangularity (\cite[Proposition 5.1]{HL15}, \cite[Theorem 3.7]{FHOO22}) between $\bbE(\bfc)$ and $ \chi_{t,q}( L(\bfc') )$ says that 
$$
\chi_{t,q}( L(\bfc) ) =  \bbE(\bfc) + \sum_{\bfc' } g_{\bfc, \bfc'}(q)
\bbE(\bfc') \qquad \text{ for some $ g_{\bfc, \bfc'}(q) \in q\Z[q]$}.
$$
Hence we have
$$
q^{N(\wt(\bfc))/2}\Omega_\calQ( \chi_{t,q}( L(\bfc) ))
\equiv  q^{N(\wt(\bfc))/2} \Omega_\calQ  (\bbE(\bfc) ) \equiv \rmG(\bfc)\ \mod q\LuphA.
$$
Since $q^{N(\wt(\bfc))/2}\Omega_\calQ( \chi_{t,q}( L(\bfc) ))$ is $c$-invariant,
we obtain
$$q^{N(\wt(\bfc))/2}\Omega_\calQ( \chi_{t,q}( L(\bfc) ))=\rmG(\bfc)$$
by Theorem~\ref{Thm: global basis}, which implies the desired result.
\end{proof}

\providecommand{\bysame}{\leavevmode\hbox to3em{\hrulefill}\thinspace}
\providecommand{\MR}{\relax\ifhmode\unskip\space\fi MR }
\providecommand{\MRhref}[2]{%
  \href{http://www.ams.org/mathscinet-getitem?mr=#1}{#2}
}
\providecommand{\href}[2]{#2}

\end{document}